\documentclass[12pt,reqno]{amsart}

\usepackage{amssymb, xcolor}
\usepackage{amscd}
\usepackage{amsfonts}
\usepackage{setspace}
\usepackage{version}
\usepackage{tikz}
\usepackage{caption}
\usepackage{hyperref}
\usepackage{comment}
\usepackage{dsfont}

\usetikzlibrary{arrows.meta}

\definecolor{regionA}{HTML}{F2BD81} 
\definecolor{regionB}{HTML}{9FCBE8} 
\definecolor{regionC}{HTML}{C6B5E2} 
\definecolor{regionD}{HTML}{A8D6AD} 
\definecolor{outline}{HTML}{1A2430}

\definecolor{squareInk}{HTML}{202C38}
\definecolor{squareMuted}{HTML}{526271}
\definecolor{squareOuter}{HTML}{44536B}
\definecolor{squareMiddle}{HTML}{006963}
\definecolor{squareInner}{HTML}{914A19}
\definecolor{squareOuterFill}{HTML}{F1F3F7}
\definecolor{squareMiddleFill}{HTML}{E5F2EF}
\definecolor{squareInnerFill}{HTML}{FFF1DE}
\definecolor{squareDot}{HTML}{8796A0}

\definecolor{ink}{HTML}{192D3C}
\definecolor{cornertext}{HTML}{526572}
\definecolor{regionone}{HTML}{F3CB70}
\definecolor{regiontwo}{HTML}{EFA384}
\definecolor{regionthree}{HTML}{E8A9C6}
\definecolor{regionfour}{HTML}{BFB5E7}
\definecolor{regionfive}{HTML}{9BCBE8}
\definecolor{regionsix}{HTML}{9CD5BD}
\definecolor{regionC}{HTML}{E6E9ED}

\definecolor{ink}{HTML}{223547}
\definecolor{circlegray}{HTML}{A4AFBB}
\definecolor{axisgray}{HTML}{CBD3DC}
\definecolor{arcblue}{HTML}{1975B8}
\definecolor{boundaryorange}{HTML}{C47A21}

\definecolor{Efill}{HTML}{CFE7F8}
\definecolor{Eink}{HTML}{16618E}
\definecolor{TOneFill}{HTML}{FFD99A}
\definecolor{TOneInk}{HTML}{A95604}
\definecolor{TTwoFill}{HTML}{C8EAD4}
\definecolor{TTwoInk}{HTML}{287C51}
\definecolor{TThreeFill}{HTML}{E2D4F3}
\definecolor{TThreeInk}{HTML}{7450A1}
\definecolor{TFourFill}{HTML}{F6CBD1}
\definecolor{TFourInk}{HTML}{B34363}
\definecolor{GapFill}{HTML}{F0F2F5}
\definecolor{Boundary}{HTML}{283847}
\definecolor{AxisColor}{HTML}{B0BAC4}
\definecolor{SquareColor}{HTML}{7B8794}

\definecolor{NEblue}{HTML}{2478B6}
\definecolor{NWgreen}{HTML}{23866E}
\definecolor{SWpurple}{HTML}{8561A8}
\definecolor{SEorange}{HTML}{CC8025}
\definecolor{InnerFill}{HTML}{F1F4F7}
\definecolor{Boundary}{HTML}{758493}
\definecolor{OuterBoundary}{HTML}{9AA6B2}
\definecolor{AxisColor}{HTML}{C1CAD3}
\definecolor{LabelColor}{HTML}{2E4052}

\usepackage{graphicx}
\usepackage{float}
\usepackage{subcaption}

\newtheorem{theorem}{Theorem}[section]
\newtheorem{proposition}[theorem]{Proposition}

\newtheorem{lemma}[theorem]{Lemma}
\newtheorem{corollary}[theorem]{Corollary}

\newtheorem{question}[theorem]{Question}

\theoremstyle{definition}

\newtheorem{example}[theorem]{Example}

\newcommand{\Ree}{\operatorname{Re}}

\hypersetup{
    colorlinks=true,
    linkcolor=blue,
    filecolor=magenta,      
    urlcolor=blue,
    pdftitle={Overleaf Example},
    pdfpagemode=FullScreen,
    }

\theoremstyle{definition}
\renewcommand{\leq}{\leqslant}
\renewcommand{\geq}{\geqslant}

\newcommand{\un}{\underline{\delta}^{-}}
\newcommand{\ov}{\overline{\delta}^{-}}

\def\R{\mathbb{R}}

\def\Z{\mathbb{Z}}

\def\T{\mathbb{T}}

\numberwithin{equation}{section}
\usepackage{amsmath, amssymb, amsthm, mathtools}
\usepackage[a4paper,margin=1in]{geometry}

\usepackage{color}
\usepackage{amssymb}

\begin{document}

\title[Discrete harmonic functions with non-negativity on a large portion of $\mathbb{Z}^2$]{On discrete harmonic functions with non-negativity on a large portion of $\mathbb{Z}^2$}
\author[J. Gaitan]{Jos\'e Gaitan Montejo}
\address[JG]{Department of Mathematics, Virginia Polytechnic Institute and State University,  225 Stanger Street, Blacksburg, VA 24061-1026, USA
}
\email{jogaitan@vt.edu}

\author[C. Garz\'on]{Carlos Garz\'on}
\address[CG]{Department of Mathematics, Virginia Polytechnic Institute and State University,  225 Stanger Street, Blacksburg, VA 24061-1026, USA
}
\email{cgarzongu@vt.edu}

\author[J. Madrid]{Jos\'e Madrid}
\address[JM]{Department of Mathematics, Virginia Polytechnic Institute and State University,  225 Stanger Street, Blacksburg, VA 24061-1026, USA
}
\email{josemadrid@vt.edu}


\date{}

\begin{abstract}
    An improvement of Liouville's theorem for harmonic functions on $\mathbb{Z}^2$ was obtained by Buhovsky, Logunov, Malinnikova, and Sodin \cite{BuhovskyLogunovMalinnikovaSodin2022}. Namely, they proved that if a discrete harmonic function is bounded on a large portion of $\mathbb{Z}^2$, then the function is constant. In this paper, we prove the corresponding result if we replace the boundedness assumption with non-negativity. 
\end{abstract}

\maketitle

\section{Introduction}

A function $u : \mathbb{Z}^2 \to \mathbb{R}$ is said to be \textit{harmonic} if, for every $x \in \mathbb{Z}^2$,
\begin{equation*}
    u(x) = \frac{1}{4} \sum_{|y|=1}u(x+y).
\end{equation*}
Additionally, for each $N \in \mathbb{Z}^+$, we consider the square 
\begin{equation*}
    Q_N = \{(n,m) \in \mathbb{Z}^2 : |n|,|m| \leq N\}.
\end{equation*}
Notice that every square $Q_N$ contains exactly $(2N+1)^2$ lattice points. For each $\varepsilon \in (0,1)$, we say that $u$ is bounded by $1$ on a $1-\varepsilon$ portion of $\mathbb{Z}^2$ if there exists $N_0 \in \mathbb{Z}^+$ such that
\begin{equation*}
    |Q_N\cap \{|u|\leq 1\}| \geq (1-\varepsilon)|Q_N|
\end{equation*}
for all $N\geq N_0$.

A classical result states that Liouville's theorem also holds in the discrete setting \cite[Theorem 5]{Heilbronn1949}. That is, if $u : \mathbb{Z}^d \to \mathbb{R}$ is a harmonic function and $|u(x)|\leq 1$ for each $x \in \mathbb{Z}^d$, then $u$ is constant. For the particular two-dimensional case, Buhovsky, Logunov, Malinnikova, and Sodin \cite{BuhovskyLogunovMalinnikovaSodin2022} obtained an unexpected improvement of this theorem. Namely, they showed the existence of a small constant $\varepsilon>0$ such that if $u : \mathbb{Z}^2 \to \mathbb{R}$ is harmonic and bounded by $1$ on a $1-\varepsilon$ portion of $\mathbb{Z}^2$, then $u$ is constant. This interesting result can be extended to periodic planar graphs \cite{Bou-RabeeCoopermanGanguly2025}, and it does not hold if we consider higher discrete dimensions or the continuous setting. Moreover, this theorem and its proof have been shown to be useful in some applications. For instance, an adaptation of its proof idea was used by Ding and Smart \cite{DingSmart2020} to establish a random quantitative unique continuation principle for eigenfunctions of Hamiltonians, which allowed them to deduce Anderson localization  for certain random Schr\"odinger operators on $\ell^2(\mathbb{Z}^2)$. Similar unique continuation results may be found in \cite{Hurtado2026,Li2022,LiZhang2022} and references therein.

Buhovsky, Logunov, Malinnikova, and Sodin's argument is based on  lower and upper exponential estimates for $u$. The lower bound \cite[Theorem B]{BuhovskyLogunovMalinnikovaSodin2022} has a continuous analogue, and its proof is deduced as a consequence  of the discrete version of Hadamard's three-circle theorem. The upper bound \cite[Theorem A]{BuhovskyLogunovMalinnikovaSodin2022} is more delicate, and it has no continuous analogue. 

Buhovsky, Logunov, Malinnikova, and Sodin also posed some open questions for discrete harmonic functions on $\mathbb{Z}^2$. One of those questions is whether the result remains true if we assume positivity rather than boundedness on a large portion of $\mathbb{Z}^2$.

\begin{question}\label{QuestionPositivity}
    Does there exist a constant $\varepsilon>0$ such that if $u : \mathbb{Z}^2 \to \mathbb{R}$ is a harmonic function that is positive on a $1-\varepsilon$ portion of $\mathbb{Z}^2$, then $u$ is constant?
\end{question}

This question motivates the results that we present in this manuscript. Before properly stating them, we introduce some additional background. For each $N \in \mathbb{Z}^+$, let 
\begin{equation*}
    \delta_N^-(u) = \frac{|\{x \in Q_N : u(x)<0\}|}{|Q_N|},
\end{equation*}
and let 
\begin{equation*}
    \underline{\delta}^-(u) = \liminf_{N\to\infty}\delta_N^-(u) \quad \text{ and } \quad \overline{\delta}^-(u) = \limsup_{N\to\infty}\delta_N^-(u).
\end{equation*}
Observe that the condition that $u$ is non-negative on a $1-\varepsilon$ portion of $\mathbb{Z}^2$ implies that $\overline{\delta}^-(u) \leq \varepsilon$. Our first main result gives a positive answer to Question \ref{QuestionPositivity} if we also assume that $u$ has at most polynomial growth. We state this result precisely in the following theorem. 

\begin{theorem}\label{thm:polynomial-main}
Let $u:\Z^2\to\R$ be a non-constant discrete harmonic function of polynomial growth.  Then, $\lim_{N\to\infty}\delta_N^-(u)$ exists and satisfies
\begin{equation}\label{eq:poly-bound-intro}
 \lim_{N\to\infty}\delta_N^-(u)\ge\sqrt2-1.
\end{equation}
Moreover, equality in \eqref{eq:poly-bound-intro} is attained by the function
\begin{equation}\label{eq:quartic-example-intro}
 U(n,m)=-(n^4-6n^2m^2+m^4)+(n^2+m^2).
\end{equation}
\end{theorem}

As a straightforward consequence of this theorem, we have the next corollary.

\begin{corollary}
Let $\varepsilon>0$. If $u :\mathbb{Z}^2 \to \mathbb{R}$ is a harmonic function with polynomial growth that is non-negative on a $2-\sqrt{2}+\varepsilon$ portion of $\mathbb{Z}^2$, then $u$ is constant.
\end{corollary}

The proof of Theorem \ref{thm:polynomial-main} essentially relies on the fact that a discrete harmonic function with polynomial growth is a polynomial itself (see \cite[Theorem 6]{Heilbronn1949}). Afterwards, we reduce the problem to proving \eqref{eq:poly-bound-intro} for the highest homogeneous part of the polynomial, and finally we use a geometrical argument in the two-dimensional setting to establish the desired inequality for homogeneous polynomials. This idea will be fully developed in Section \ref{SectionPolynomial}.

By using Theorem \ref{thm:polynomial-main}, we can prove the main result of this manuscript, which gives a positive answer to Question \ref{QuestionPositivity} in its full generality. 

\begin{theorem}\label{thm:QuestionAnswer}
    There exists a sufficiently small constant $\varepsilon_{\ast}>0$ such that every harmonic function $u : \mathbb{Z}^2 \to \mathbb{R}$ with $\overline{\delta}^-(u) \leq \varepsilon_{\ast}$ is constant. In particular, if $u$ is a harmonic function on $\mathbb{Z}^2$ and is non-negative on a $1-\varepsilon_{\ast}$ portion of $\mathbb{Z}^2$, then $u$ is constant. 
\end{theorem}

Our approach to deduce this result is similar to that of \cite[Theorem 1.1]{BuhovskyLogunovMalinnikovaSodin2022}. Namely, we first find an exponential lower bound for the growth of a non-constant harmonic function on $\mathbb{Z}^2$ for which $\overline{\delta}^-(u)$ is sufficiently small.

\begin{theorem}\label{thm:LowerExponentialBound}
    There exist constants $\delta_L,b > 0$ such that if $u : \mathbb{Z}^2 \to \mathbb{R}$ is a non-constant harmonic function and $\overline{\delta}^-(u) \leq \delta_L$, then there exists $R_0 = R_0(u) \in \mathbb{Z}^+$ satisfying
    \begin{equation}\label{eq:LowerExpGrowthU}
        \|u\|_{Q_R} \geq e^{bR}
    \end{equation}
    for every $R\geq R_0$.
\end{theorem}

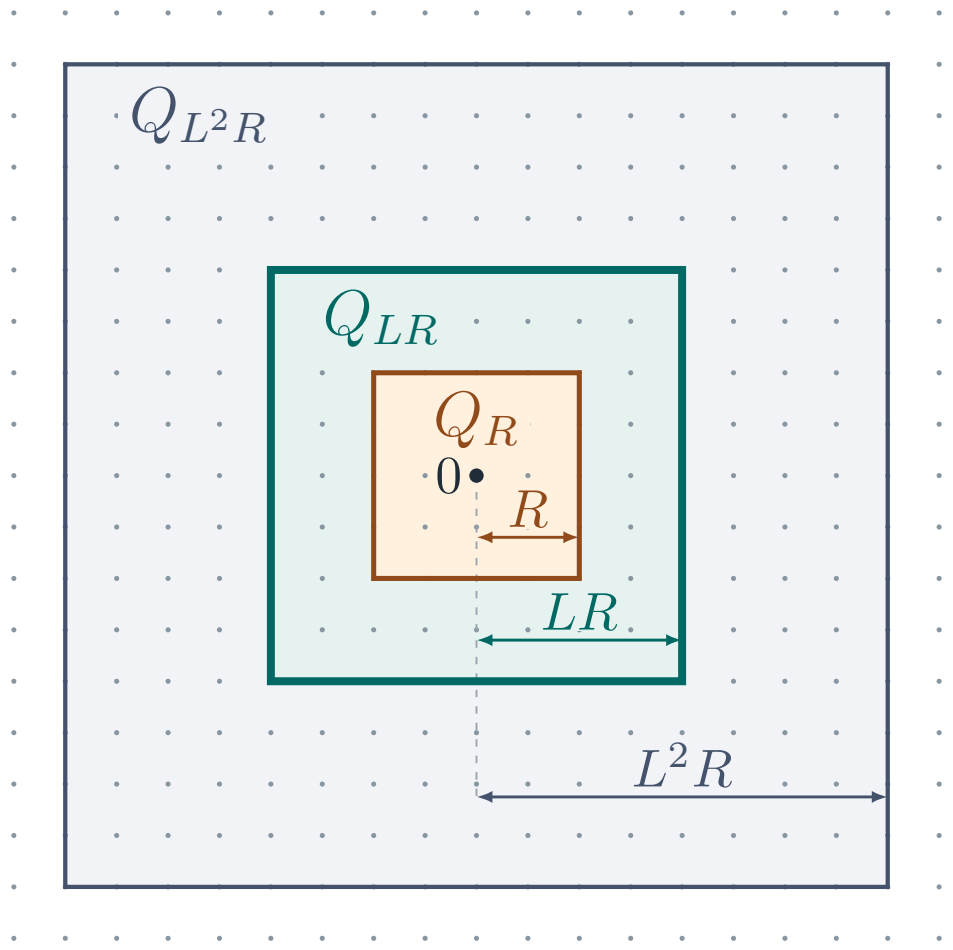
\begin{figure}[H]
    \centering
    \begin{tikzpicture}[scale=1,x=0.34cm, y=0.34cm, text=squareInk, line join=miter,
  boundary label/.style={inner xsep=2pt, inner ysep=0.6pt, font=\normalsize},
  distance label/.style={midway, above=1.4pt, inner xsep=1pt,
    inner ysep=0.3pt, font=\footnotesize},
  dimension line/.style={<->, >= {Latex[length=3pt,width=2.3pt]},
    line width=0.55pt}]

  \fill[squareOuterFill]  (-8,-8) rectangle (8,8);
  \fill[squareMiddleFill] (-4,-4) rectangle (4,4);
  \fill[squareInnerFill]  (-2,-2) rectangle (2,2);
  \foreach \i in {-9,...,9}{
    \foreach \j in {-9,...,9}{
      \fill[squareDot] (\i,\j) circle[radius=0.48pt];
    }
  }
  \draw[squareMuted!58,line width=0.35pt,dash pattern=on 1.4pt off 1.7pt]
    (0,0) -- (0,-6.25);
  \draw[squareOuter, line width=0.8pt] (-8,-8) rectangle (8,8);
  \draw[squareMiddle,line width=1.5pt] (-4,-4) rectangle (4,4);
  \draw[squareInner, line width=0.95pt] (-2,-2) rectangle (2,2);

  \node[boundary label,fill=squareOuterFill,text=squareOuter]
    at (-5.4,7) {$Q_{L^2R}$};
  \node[boundary label,fill=squareMiddleFill,text=squareMiddle]
    at (-1.85,3.05) {$Q_{LR}$};
  \node[boundary label,fill=squareInnerFill,text=squareInner]
    at (0,1.08) {$Q_R$};

  \draw[dimension line,squareInner] (0,-1.20)
    -- node[distance label,fill=squareInnerFill] {$R$} (2,-1.20);
  \draw[dimension line,squareMiddle] (0,-3.20)
    -- node[distance label,fill=squareMiddleFill] {$LR$} (4,-3.20);
  \draw[dimension line,squareOuter] (0,-6.25)
    -- node[distance label,fill=squareOuterFill] {$L^2R$} (8,-6.25);
  \fill[squareInk] (0,0) circle[radius=1.35pt];
  \node[anchor=east,inner sep=0.6pt,font=\footnotesize,fill=squareInnerFill]
    at (-0.20,0) {$0$};
\end{tikzpicture}
    \caption{Discrete three-circle theorem.}
    \label{fig:ThreeSquares}
\end{figure}

This result is analogous to \cite[Theorem B]{BuhovskyLogunovMalinnikovaSodin2022}. In fact, its proof also relies on establishing a suitable discrete version of Hadamard's three circle theorem, which allows us to estimate $\|u\|_{Q_{LR}}$ by a combination of an exponential decay multiple of $\|u\|_{Q_{L^2 R}}$ and an exponential growth multiple of $\|u\|_{Q_{R}}$ (see Figure \ref{fig:ThreeSquares} for an illustration of these squares). The argument is explained in full detail in Section \ref{Section:ExpLowerBound}. Hence, our goal now is to obtain an analogous result to the exponential upper bound \cite[Theorem A]{BuhovskyLogunovMalinnikovaSodin2022} for harmonic functions whose negative set has sufficiently small upper density on $\mathbb{Z}^2$. For $\overline{\delta}^-(u)<1/32$, the following proposition gives a partial exponential upper bound result.

\begin{proposition}\label{prop:large-range}
There is an absolute constant $C<\infty$ independent of $u$ such that every discrete harmonic function $u:\Z^2\to\R$ with $\delta=\ov(u)<1/32$ satisfies
\begin{equation}\label{eq:main-bound}
 |u(x,y)|\le C_u\exp\!\left(C\log\frac{e}{1-32\delta}\,(|x|+|y|)\right),
\end{equation}
for some $C_u<\infty$, and for all $(x,y)\in\Z^2$. 
\end{proposition}

Proposition~\ref{prop:large-range} is only a partial analogue to \cite[Theorem A]{BuhovskyLogunovMalinnikovaSodin2022} because the exponential rate does not decrease to zero as $\delta\to0$. However, for a carefully chosen collection $G$ of centers such that the negative values of $u$ occupy at most a sufficiently small fixed fraction of every square centered at $z\in G$, it is possible to obtain estimates such as \eqref{eq:main-bound} in translated diamonds $D_R(z)=\{x\in\Z^2:|x-z|_1\le R\}$ with bounds that grow at the same rate as the magnitude of $u$ at the center $z\in G$. Taking the infimum of these bounds over all such centers and controlling the size of the complement $\Z^2\setminus G$ (using Vitali's covering lemma) lead us to our main exponential upper bound result. In the next theorem, we obtain an improvement of Proposition~\ref{prop:large-range} that gives a complete analogue of \cite[Theorem A]{BuhovskyLogunovMalinnikovaSodin2022} under the non-negativity condition. 

\begin{theorem}\label{thm:main}
There are absolute constants $0<\delta_U<1/32$ and $C<\infty$ with the following property: If $u:\Z^2\to\R$ is discrete harmonic and $\overline{\delta}^-(u) \leq \delta_U$, there exist constants $C_u,R_u$ such that
\begin{equation}\label{eq:sparse-finite}
 \log(1+\|u\|_{Q_R})\leq C_u+CR\bigl(\delta_{6R}^-(u)\bigr)^{1/3}
\end{equation}
for all $R\geq R_u$. In particular,
\begin{equation}\label{eq:sparse-rate}
 \limsup_{R\to\infty}\frac{\log(1+\|u\|_{Q_R})}{R}\leq C\overline{\delta}^-(u)^{1/3}.
\end{equation}
\end{theorem}

Observe that the exponent in \eqref{eq:sparse-rate} tends to zero with the negative density $\overline{\delta}^-(u)$. In particular, $\overline{\delta}^-(u)=0$ implies that $\log(1+\|u\|_{Q_R})=o(R)$. The proof of Proposition~\ref{prop:large-range} uses an extension argument that controls the maxima of $|u|$ on successively larger diamonds centered at the origin via a growth estimate on an intermediate square, see Figure \ref{fig1}. We show a rough sketch of this proof, where we point out the main steps and reductions in the argument. The proof is given in full detail in Section \ref{Section:ExpGrowth}. Before the outline of the proof, we show how Theorem \ref{thm:QuestionAnswer} follows as a straightforward consequence of Theorems \ref{thm:LowerExponentialBound} and \ref{thm:main}.

\begin{proof}[Proof of Theorem \ref{thm:QuestionAnswer}]
    Let $C$ and $b$ be the constants from \eqref{eq:sparse-rate} and \eqref{eq:LowerExpGrowthU}, respectively. Let 
    \begin{equation*}
        \varepsilon_{\ast} = \min\left\{\delta_L , \delta_U , \left(  \frac{b}{2C}\right)^3\right\} > 0.
    \end{equation*}
    Let $u$ be a harmonic function on $\mathbb{Z}^2$ such that $\overline{\delta}^-(u) \leq \varepsilon_{\ast}$. For the sake of contradiction, let us assume that $u$ is non-constant. By Theorems \ref{thm:LowerExponentialBound} and \ref{thm:main}, we deduce
    \begin{equation*}
        b \leq \liminf_{R \to \infty} \frac{\log(1+\|u\|_{Q_R})}{R} \leq \limsup_{R\to\infty}\frac{\log(1+\|u\|_{Q_R})}{R}\leq C\overline{\delta}^-(u)^{1/3} \leq C \varepsilon_{\ast}^{1/3}.
    \end{equation*}
    However, this implies that $\varepsilon_{\ast} \geq (b/C)^3$, which contradicts the choice of $\varepsilon_{\ast}$. Therefore, $u$ is constant.
\end{proof}

\subsection{Outline of the proof of Proposition~\ref{prop:large-range}}

\begin{figure}[H]
    \centering
    \begin{tikzpicture}[
    x=0.4cm, y=0.4cm,
    scale=0.6,
    transform shape,
    line cap=round,
    line join=round,
    >={Stealth[length=4pt]},
    leader/.style={
        line width=0.6pt,
        shorten <=3pt,
        shorten >=3pt
    }
]

\def\R{10}
\def\K{4}
\pgfmathtruncatemacro{\outerR}{\R+\K}
\pgfmathtruncatemacro{\squareR}{\R-1}

\fill[InnerFill]
    (\R,0) -- (0,\R) --
    (-\R,0) -- (0,-\R) -- cycle;

\foreach \sx/\sy/\regioncolor in {
    1/1/NEblue,
    -1/1/NWgreen,
    -1/-1/SWpurple,
    1/-1/SEorange%
}{
    \begin{scope}[xscale=\sx,yscale=\sy]
        \fill[\regioncolor,opacity=0.055]
            (1,\squareR)
            -- ({\K+1},\squareR)
            -- (\squareR,{\K+1})
            -- (\squareR,1) -- cycle;
    \end{scope}
}

\draw[->,AxisColor,line width=0.5pt]
    ({-\outerR-1.7},0) -- ({\outerR+1.9},0);
\draw[->,AxisColor,line width=0.5pt]
    (0,{-\outerR-1.7}) -- (0,{\outerR+1.9});

\node[Boundary,anchor=west]
    at ({\outerR+2.1},0.24) {\Large $x$};
\node[Boundary,anchor=south west]
    at (0.22,{\outerR+2}) {\Large $y$};

\draw[
    Boundary,
    line width=0.7pt,
    dash pattern=on 4pt off 3pt
]
    (-\squareR,-\squareR)
    rectangle (\squareR,\squareR);

\draw[Boundary,line width=0.8pt]
    (\R,0) -- (0,\R) --
    (-\R,0) -- (0,-\R) -- cycle;

\draw[OuterBoundary,line width=0.9pt]
    (\outerR,0) -- (0,\outerR) --
    (-\outerR,0) -- (0,-\outerR) -- cycle;

\foreach \sx/\sy/\regioncolor in {
    1/1/NEblue,
    -1/1/NWgreen,
    -1/-1/SWpurple,
    1/-1/SEorange%
}{
    \begin{scope}[xscale=\sx,yscale=\sy]
        \foreach \kk in {1,...,\K}{
            \pgfmathtruncatemacro{\firstJ}{\kk+1}

            \draw[\regioncolor,line width=1.2pt]
                (\firstJ,\squareR)
                -- (\squareR,\firstJ);

            \foreach \jj in {\firstJ,...,\squareR}{
                \filldraw[
                    fill=\regioncolor,
                    draw=white,
                    line width=0.2pt
                ]
                    (\jj,{\R+\kk-\jj})
                    circle[radius=1.35pt];
            }
        }
    \end{scope}
}

\node[
    Boundary,
    font=\LARGE,
    fill=InnerFill,
    inner sep=1pt
] at (-1.8,1.1) {$D_R$};

\node[Boundary,font=\small,anchor=north east]
    at (-0.25,-0.38) {$0$};

\node[LabelColor,font=\Large,anchor=south] (Dlabel)
    at ({-\outerR+1.6},{\outerR-1.9}) {$D_{R+K}$};

\draw[leader,LabelColor]
    (Dlabel.south east)
    -- ({-0.54*\outerR},{0.46*\outerR});

\node[Boundary,font=\large,anchor=south] (Qlabel)
    at ({\squareR+2.7},{\squareR+2.4}) {\Large $Q_{R-1}$};

\draw[leader,Boundary]
    (Qlabel.south west) -- (\squareR,\squareR);

\end{tikzpicture}
    \caption{Illustration of $(D_{R+K}\cap Q_{R-1}) \setminus D_R$.}
    \label{fig1}
\end{figure}
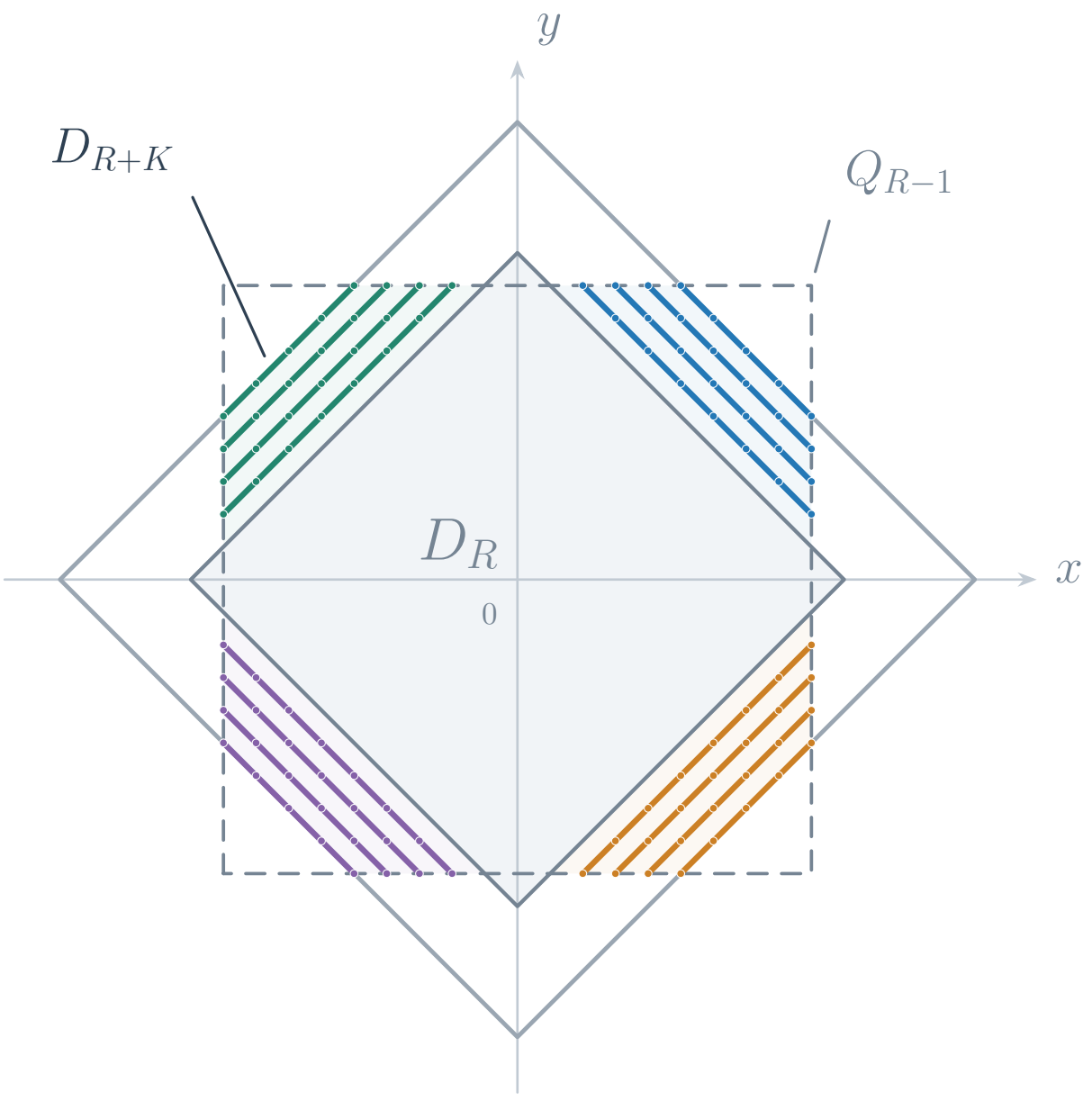

For each $R \in \mathbb{Z}^+$, let 
\begin{equation*}
    D_R = \{(x,y) \in \mathbb{Z}^2 : |x| + |y| \leq R\} \quad \text{and} \quad M_u(D_R) = \max\left\{1,\|u\|_{D_R}\right\}.
\end{equation*}
Moreover, let $\gamma = (1-32\delta)/8$, where $\delta := \overline{\delta}^{-}(u) < 1/32$. Because of the assumption on $\delta$, we have that $\gamma \in (0,1/8]$. Our initial step is to find a suitable reduction for \eqref{eq:main-bound}. Indeed, we show that if for every $R \in \mathbb{Z}^+$ large enough there exists $K \in [\gamma R , (1/2-\gamma)R] \cap \mathbb{Z}$ for which
\begin{equation*}
     \|u\|_{D_{R+K}} \leq M_u(D_R) e^{C K \log(eR/K)}
\end{equation*}
holds, then we can iterate this inequality to deduce \eqref{eq:main-bound}. By using the harmonicity of $u$, we can prove that
\begin{equation*}
    M_u(D_{R+K}) \leq 7^{K+1} M_u(D_{R+K} \cap Q_{R-1}).
\end{equation*}
Observe that $7^{K+1} \leq e^{2K \log 7} \leq  e^{\widetilde{C} K \log(eR/K)}$. Hence, it remains to estimate the maximum of $|u|$ in $D_{R+K} \cap Q_{R-1}$ by an exponential multiple of $M_u(D_R)$. In other words, we need to estimate $|u|$ on the lines composing the set $(D_{R+K}\cap Q_{R-1}) \setminus D_R$, which are illustrated in Figure \ref{fig1}.

By means of a rotation argument, it is enough to show the desired estimate for the highlighted lines from Figure \ref{fig1} that lie in the first quadrant. Applying the harmonicity of $u$ once again, we can find a suitable decomposition of $u$ on each of these lines as a sum of two terms, where the first one can be easily bounded in terms of $M_u(D_R)$ and the second one is a polynomial. The harmonicity of $u$ yields a convenient relation between the polynomials corresponding to three consecutive lines, and thus it suffices to estimate the polynomial associated to the line $\ell_K =\{x+y = R+K\} \cap Q_{R-1}$. This is the point where the non-negativity hypothesis plays an essential role, since we prove that the desired bound for such a polynomial follows from ensuring that there are at most $R/2 - K - 1$ points in $\ell_K$ where $u$ takes negative values. By the non-negativity assumption, for each $R$ large enough, we can choose $K \in [\gamma R , (1/2-\gamma)R] \cap \mathbb{Z}$ satisfying this condition for the last line in each quadrant, which finishes the proof.

\subsection{Organization of the paper} In Section \ref{SectionPolynomial}, we prove Theorem \ref{thm:polynomial-main}, which provides a positive answer to Question \ref{QuestionPositivity} for functions with polynomial growth. In Section \ref{Section:ExpLowerBound}, we show Theorem \ref{thm:LowerExponentialBound}, which gives an exponential lower bound for non-constant discrete harmonic functions which are non-negative on a large portion of $\mathbb{Z}^2$. Finally, in Section  \ref{Section:ExpGrowth}, we prove an exponential upper bound for the growth of such functions, which is stated in Proposition \ref{prop:large-range}, and then we use this result to obtain its stronger version contained in Theorem \ref{thm:main}. 

\section{Polynomial growth: Proof of Theorem \ref{thm:polynomial-main}} \label{SectionPolynomial}

A classical result  \cite[Theorem 6]{Heilbronn1949} states that a discrete harmonic function of polynomial growth is a polynomial. Hence, to prove Theorem \ref{thm:polynomial-main}, we begin by assuming that $u$ is a polynomial. In the next lemma, we use some elementary facts about harmonic polynomials to obtain a suitable representation for the top degree homogeneous part of $u$.

\begin{lemma}\label{lem:leading-harmonic}
Let $u$ be a nonconstant discrete harmonic polynomial on $\mathbb{Z}^2$ of degree $d$, and let $H_d$ be its highest homogeneous part.  Then, $H_d$ is a nonzero homogeneous Euclidean harmonic polynomial. Consequently,
\begin{equation}\label{eq:harmonic-homogeneous-form}
 H_d(r\cos\theta,r\sin\theta)
 =A r^d\cos(d\theta-\varphi)
\end{equation}
for some $A\ne0$ and $\varphi\in\R$.
\end{lemma}
\begin{proof}
    It is well known that, for any harmonic polynomial in $\mathbb{Z}^n$, its top degree homogeneous part  is a harmonic polynomial in $\mathbb{R}^n$ (see \cite{LippnerMangoubi2015}). Moreover, a homogeneous harmonic polynomial of degree $d$ restricted to the unit circle belongs to $\text{span}\{\cos(d\theta),\sin(d\theta)\}$ (see \cite[p. 82]{AxlerBourdonRameyBook}). It follows that $H_d(\cos\theta,\sin\theta) = A \cos(d\theta-\varphi)$ for some $A \in \mathbb{R}\setminus \{0\}$, which, combined with the homogeneity of $H_d$, immediately implies \eqref{eq:harmonic-homogeneous-form}.
\end{proof}

Next, we relate the asymptotic behavior of $\delta_N^-(u)$ to the area of the subset of $[-1,1]^2$ where the highest homogeneous term of $u$ takes negative values. As a result, the problem is reduced to estimating the area of such a set. 

\begin{lemma}\label{lemma:LimitDiscretePols}
Let $u$ be a nonconstant discrete harmonic polynomial on $\mathbb{Z}^2$, and let $H_d$ be its highest homogeneous part. Then, $\lim_{N\to \infty} \delta_N^-(u)$ exists and
\begin{equation}\label{eq:ConvDeltN-}
    \lim_{N\to \infty} \delta_N^-(u) = \frac{|\{H_d < 0\} \cap [-1,1]^2|}{4}.
\end{equation}
\end{lemma}
\begin{proof}
    We start by proving that 
    \begin{equation*}
        N^{-d}u(\lfloor Nx \rfloor , \lfloor Ny \rfloor) \to H_d(x,y)
    \end{equation*}
    uniformly for $(x,y) \in [-1,1]^2$.

    Notice that $u$ can be decomposed as
    \begin{equation*}
        u(x,y) = \sum_{j=0}^d H_j(x,y),
    \end{equation*}
    where each $H_j$ is a homogeneous polynomial of degree $j$. We can choose a constant $L$ sufficiently large that 
    \begin{equation*}
        \sup_{0\leq j\leq d}\{\|H_j\|_{L^{\infty}([-1,1]^2)},\|\nabla H_j\|_{L^{\infty}([-1,1]^2)}\} \leq L.
    \end{equation*}
    Let $(x,y) \in[-1,1]^2$ and $N \in \mathbb{N}$, and define $x_N := \lfloor Nx \rfloor/N$ and $y_N := \lfloor Ny \rfloor/N$. Notice that $(x_N,y_N) \in [-1,1]^2$. Because of the homogeneity of each $H_j$, an application of the triangle inequality yields
    \begin{equation*}
        \begin{aligned}
            |N^{-d}u(\lfloor Nx \rfloor , \lfloor Ny \rfloor) - H_d(x,y)| &= |N^{-d} u(Nx_N , Ny_N) - H_d(x,y)|\\ &\leq |H_d(x_N,y_N) - H_d(x,y)| + L \sum_{j=0}^{d-1} \frac{1}{N^{d-j}}.
        \end{aligned}
    \end{equation*}
     Since $\lfloor Nx \rfloor \leq Nx \leq \lfloor Nx \rfloor + 1$, we have that $x_N \leq x \leq x_N + 1/N$. This implies that
    \begin{equation*}
        |x_N - x| = x - x_N \leq \frac{1}{N}.
    \end{equation*}
    The same argument shows that $|y_N - y| \leq 1/N$. As a result,
    \begin{equation*}
        \|(x_N,y_N) - (x,y)\| \leq \frac{\sqrt{2}}{N}.
    \end{equation*}
    By the mean value estimate, we deduce
    \begin{equation*}
        |N^{-d}u(\lfloor Nx \rfloor , \lfloor Ny \rfloor) - H_d(x,y)| \leq L \left( \frac{\sqrt{2}}{N} + \sum_{j=0}^{d-1} \frac{1}{N^{d-j}} \right).
    \end{equation*}
    Therefore,
    \begin{equation*}
        \|N^{-d}u(\lfloor N (\cdot) \rfloor , \lfloor N(\cdot )\rfloor) - H_d\|_{L^{\infty}([-1,1]^2)} \to 0
    \end{equation*}
    as $N\to \infty$. This shows the desired uniform convergence.

    Now we shall use this convergence to deduce \eqref{eq:ConvDeltN-}.

    For each $N \in \mathbb{N}$, let us consider $F_N(x,y) = N^{-d}u(\lfloor Nx \rfloor , \lfloor Ny \rfloor)$ and $\varepsilon_N = \|F_N - H_d\|_{L^{\infty}([-1,1]^2)}$. Notice that if $H_d(x,y) > \varepsilon_N$, then
    \begin{equation*}
        F_N(x,y) \geq H_d(x,y) - \varepsilon_N > 0.
    \end{equation*}
    Similarly, if $H_d(x,y) < - \varepsilon_N$, then
    \begin{equation*}
        F_N(x,y) \leq H_d(x,y) + \varepsilon_N < 0.
    \end{equation*}
    This implies that the only points in the square $[-1,1]^2$ where $F_N$ and $H_d$ may have different signs are contained in $\{|H_d| \leq \varepsilon_N\} \cap [-1,1]^2$. Thus,
    \begin{equation*}
        |m(\{F_N < 0\} \cap [-1,1]^2) - m(\{H_d < 0\} \cap [-1,1]^2)| \leq m(\{|H_d| \leq \varepsilon_N\} \cap [-1,1]^2).
    \end{equation*}
    Applying the dominated convergence theorem and the fact that the zero set of $H_d$ has null measure, we deduce
    \begin{equation*}
        m(\{F_N < 0\} \cap [-1,1]^2) \to m(\{H_d < 0\} \cap [-1,1]^2)
    \end{equation*}
    as $N \to \infty$. Next, for any $i,j \in \{-N,\dots,N-1\}$, let us consider 
    \begin{equation*}
        C_{i,j} = \left[ \frac{i}{N} , \frac{i+1}{N} \right) \times \left[ \frac{j}{N} , \frac{j+1}{N} \right).
    \end{equation*}
    Observe that these squares are a partition of $[-1,1)^2$. Moreover, for each $(x,y) \in C_{i,j}$, we get that $F_N(x,y) = N^{-d}u(i,j)$. This means $F_N$ has the same sign in the whole square $C_{i,j}$ as $u(i,j)$. Because the area of each of these squares is $N^{-2}$, we get
    \begin{equation*}
        m(\{F_N < 0\} \cap [-1,1]^2) = N^{-2}\cdot |\{(i,j) \in \{-N,\dots,N-1\}^2 : u(i,j) < 0\}|.
    \end{equation*}
    As a result, for each $N$ there exists $r_N \in [0,4N+1]$ such that
    \begin{equation*}
    \begin{aligned}
        |\{u<0\}\cap Q_N| &= |\{(i,j) \in \{-N,\dots,N-1\}^2 : u(i,j) < 0\}| + r_N\\ &= N^2 \cdot m(\{F_N < 0\} \cap [-1,1]^2) + r_N.
    \end{aligned}
    \end{equation*}
    Given that $|Q_N| = (2N+1)^2$, we conclude
    \begin{equation*}
        \lim_{N\to \infty} \frac{|\{u<0\}\cap Q_N|}{|Q_N|} = \frac{m(\{H_d < 0\} \cap [-1,1]^2)}{4}.
    \end{equation*}
\end{proof}

Finally, we find sharp lower bounds for the area of the subset of $[-1,1]^2$ where a given homogeneous harmonic polynomial takes negative values. In this proof, the representation obtained in Lemma \ref{lem:leading-harmonic} plays an essential role, as well as a proper estimate for the area of the desired set by means of polar coordinates in $\mathbb{R}^2$.

\begin{lemma}\label{lem:square-area}
Let $H$ be a nonzero homogeneous harmonic polynomial of degree $d\ge1$ on $\R^2$.  Then
\begin{equation}\label{eq:square-area-bound}
 \frac{|\{H<0\}\cap[-1,1]^2|}{4}
 \ge
 \begin{cases}
 \frac12,&4\nmid d,\\[1mm]
 \sqrt2-1,&4\mid d.
 \end{cases}
\end{equation}
Moreover, the lower bound $\sqrt2-1$ is attained when $d=4$ and
\begin{equation}\label{eq:OptimalHLemma2.3}
    H(x,y)=-\Ree(x+iy)^4=-(x^4-6x^2y^2+y^4).
\end{equation}
\end{lemma}
\begin{proof}
Applying the second part of Lemma \ref{lem:leading-harmonic}, we can represent $H$ as $H(r\cos\theta,r\sin\theta) = Ar^d\cos(d\theta-\varphi)$ for some $A>0$ and $\varphi\in\mathbb{R}$.

To begin with, let us assume that $d$ is odd. Because of the homogeneity of $H$, we have that $H(-x,-y) = (-1)^d H(x,y) = -H(x,y)$. Let $P = \{H>0\}\cap [-1,1]^2$ and $N=\{H<0\}\cap [-1,1]^2$. Notice that $(x,y) \in P$ if and only if $(-x,-y) \in N$. Letting $T : \mathbb{R}^2 \to \mathbb{R}^2$ be the linear operator defined as
\begin{equation*}
    T\begin{pmatrix}
        x\\y
    \end{pmatrix} = \begin{pmatrix}
        -x\\-y
    \end{pmatrix} = \begin{pmatrix}
        \cos(\pi) & \sin(\pi)\\ -\sin(\pi) & \cos(\pi)
    \end{pmatrix} \begin{pmatrix}
        x\\ y
    \end{pmatrix},
\end{equation*}
we get that $T$ is the rotation of angle $\pi$ and $T(P) = N$. Since $T$ preserves areas and the zero set of $H$ has null area, we get that $|P| = |N|$ and $|P\cup N| = |[-1,1]^2| = 4$, which implies that $|N| = |P| = 2$. As a result,
\begin{equation*}
    \frac{|\{H<0\}\cap[-1,1]^2|}{4} = \frac{|N|}{4} = \frac{1}{2}
\end{equation*}
when $d$ is odd.

If $d$ is even but not a multiple of $4$ (i.e. $d\equiv 2 \mod 4$), then $d\pi/2$ is an odd multiple of $\pi$. As a result,
\begin{equation}\label{eq:HinvRot}
    \begin{aligned}
        H(-y,x) &= H(r\cos(\theta+\pi/2) , r\sin(\theta+\pi/2))\\ &= Ar^d \cos(d\theta -\varphi + d\pi/2)\\ &=  Ar^d \cos(d\theta -\varphi + \pi)\\ &= - Ar^d \cos(d\theta -\varphi)\\ &= -H(x,y). 
    \end{aligned}
\end{equation}
By considering the counterclockwise rotation of angle $\pi/2$ and arguing as in the previous case, we deduce that 
\begin{equation*}
    \frac{|\{H<0\}\cap[-1,1]^2|}{4} = \frac{1}{2}
\end{equation*}
when $d\equiv 2 \mod 4$.

Now let us assume that $d$ is a multiple of $4$. Observe that the square $[-1,1]^2$ can be decomposed into four congruent regions as follows:

\begin{figure}[H]
    \centering
    \begin{tikzpicture}[
    scale=2,
    text=outline,
    line join=miter,
    region label/.style={font=\sffamily\large}
]
    \fill[regionA] (0,0) -- (1,-1) -- (1,1) -- cycle;
    \fill[regionB] (0,0) -- (1,1) -- (-1,1) -- cycle;
    \fill[regionC] (0,0) -- (-1,1) -- (-1,-1) -- cycle;
    \fill[regionD] (0,0) -- (-1,-1) -- (1,-1) -- cycle;

    \draw[outline,line width=1pt]
        (-1,-1) rectangle (1,1);

    \draw[outline,line width=0.8pt]
        (-1,-1) -- (1,1)
        (-1,1) -- (1,-1);

    \fill[outline] (0,0) circle[radius=0.012];

    \node[region label] at ( 0.67, 0   ) {A};
    \node[region label] at ( 0,    0.67) {B};
    \node[region label] at (-0.67, 0   ) {C};
    \node[region label] at ( 0,   -0.67) {D};

    \node[above=5pt] at (-1, 1) {$(-1,1)$};
    \node[above=5pt] at ( 1, 1) {$(1,1)$};
    \node[below=5pt] at (-1,-1) {$(-1,-1)$};
    \node[below=5pt] at ( 1,-1) {$(1,-1)$};
\end{tikzpicture}
    \caption{Regions $A,B,C,D$ in the square $[-1,1]^2$.}
    \label{fig:placeholder1}
\end{figure}
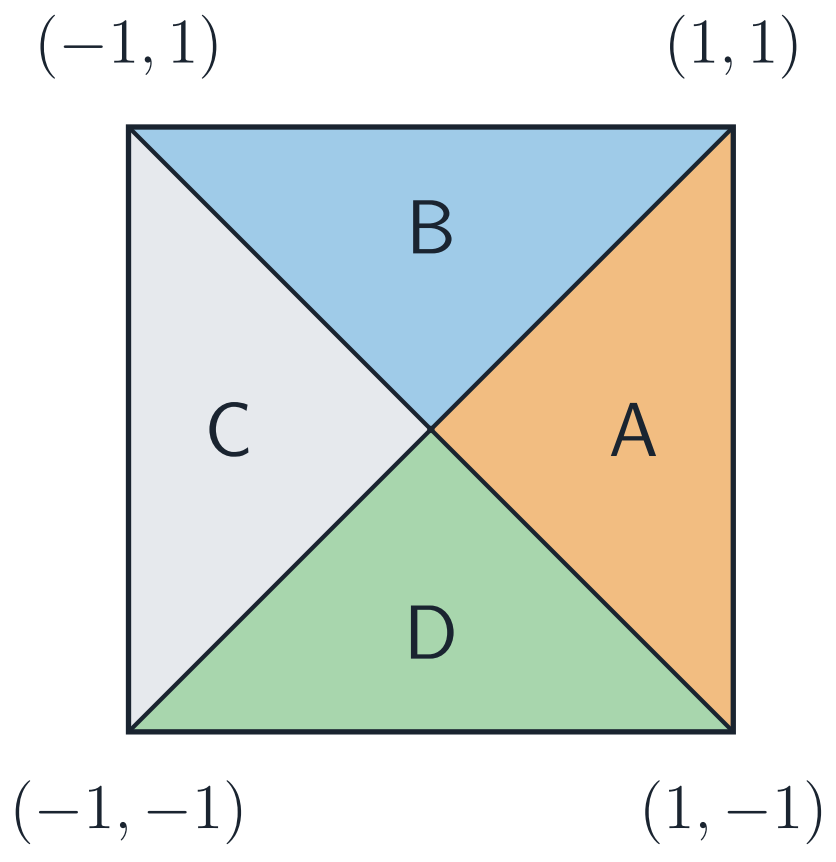
Notice that the region $A$ can be written in polar coordinates as
\begin{equation}\label{eq:PolarCoordinatesA}
    A = \{(r\cos\theta,r\sin\theta) : -\pi/4 \leq \theta \leq \pi/4 \text{ and }  0\leq r\leq \sec\theta\}.
\end{equation}
Observe that the period of $\cos(d\theta-\varphi)$ is $2\pi/d = (\pi/2)/(d/4)$, which implies that the interval $[-\pi/4,\pi/4]$ contains $d/4$ periods of $\cos(d\theta-\varphi)$. Hence, the set 
\begin{equation*}
    E = \{\theta \in [-\pi/4,\pi/4] : \cos(d\theta - \varphi) < 0\}
\end{equation*}
has length $\pi/4$. Notice that for $r>0,$ $(x,y) = (r\cos\theta,r\sin\theta) \in A \cap \{H<0\}$ if and only if $\theta \in E$. Applying the area formula in polar coordinates, we get
\begin{equation*}
    |A \cap \{H<0\}| = \frac{1}{2}\int_{E}\sec^2\theta \, d\theta.
\end{equation*}
Since $\sec^2\theta$ is even and increasing in $|\theta|$ in the interval $[-\pi/4,\pi/4]$, it follows that
\begin{equation*}
    \min\left\{\int_{E}\sec^2\theta d\theta : E \subseteq [-\pi/4,\pi/4] \text{ and } |E| = \pi/4 \right\} = \int_{-\pi/8}^{\pi/8} \sec^2\theta \, d\theta = 2\tan(\pi/8). 
\end{equation*}
Thus,
\begin{equation*}
    |A \cap \{H<0\}| \geq \tan(\pi/8) = \sqrt{2}-1.
\end{equation*}
Because $d$ is a multiple of $4$, we can argue as in \eqref{eq:HinvRot} to deduce that $H(-y,x) = H(x,y)$. This means that $H$ is invariant under counterclockwise rotation of angle $\pi/2$, which implies that 
\begin{equation*}
    |X \cap \{H<0\}| = |A \cap \{H<0\}| \geq  \sqrt{2}-1
\end{equation*}
for all regions $X \in \{B,C,D\}$. Therefore,
\begin{equation*}
    \frac{|\{H<0\}\cap[-1,1]^2|}{4} \geq \sqrt{2}-1
\end{equation*}
when $4 \mid d$. This finishes the proof of \eqref{eq:square-area-bound}.

Finally, we shall show that the left-hand side of \eqref{eq:square-area-bound} is equal to $\sqrt{2}-1$ for the function $H$ defined in \eqref{eq:OptimalHLemma2.3}, which clearly is harmonic on $\mathbb{R}^2$ since it is the real part of a holomorphic function on $\mathbb{C}$.

Observe that the rotation of angle $\pi/2$ is given by $z\mapsto iz$ in complex notation. Since $(iz)^4 = z^4$, it follows that $H$ is preserved by this rotation. Hence, it is enough to calculate the area of $\{H<0\}$ in the region $A$ of the square $[-1,1]^2$ (see Figure \ref{fig:placeholder1}). By using the polar decomposition of $A$ given in \eqref{eq:PolarCoordinatesA} and the fact that $H(r\cos\theta,r\sin\theta) = - r^4 \cos(4\theta)$, we deduce
\begin{equation*}
    \{H<0\}\cap A = \{re^{i\theta} \in A \setminus \{0\} : \cos(4\theta)>0\} = \{re^{i\theta} \in A \setminus \{0\} : -\pi/8< \theta< \pi/8\}.
\end{equation*}
By using the area formula for polar coordinates as before, we conclude
\begin{equation*}
    |\{H<0\}\cap A| = \frac{1}{2}\int_{-\pi/8}^{\pi/8} \sec^2\theta \, d\theta = \sqrt{2}-1.
\end{equation*}

\end{proof}

Observe that \eqref{eq:poly-bound-intro} follows as a straightforward consequence of Lemmas \ref{lem:leading-harmonic}, \ref{lemma:LimitDiscretePols}, and \ref{lem:square-area}. Moreover, notice that the leading homogeneous term for the function $U$ defined in \eqref{eq:quartic-example-intro} is the polynomial $H$ defined in \eqref{eq:OptimalHLemma2.3}. Thus, in order to conclude the proof of Theorem \ref{thm:polynomial-main}, it remains to see that this function $U$ is discrete harmonic. Since this can be verified by using some basic algebraic manipulations, we omit the details.

\section{Exponential lower bound: Proof of Theorem \ref{thm:LowerExponentialBound}} \label{Section:ExpLowerBound}

From Theorem \ref{thm:polynomial-main}, we know that a non-constant discrete harmonic function which is non-negative on a sufficiently large portion of $\mathbb{Z}^2$ cannot have polynomial growth. By using this idea, we can reduce the desired exponential lower bound of Theorem \ref{thm:LowerExponentialBound} to a discrete version of Hadamard's three circle theorem. This strategy can be seen as an analogue of that used in \cite[Section 2]{BuhovskyLogunovMalinnikovaSodin2022}.

\begin{lemma}\label{lemma:ReductionLower}
    Let $L \in \mathbb{Z}$ and $\tau,\kappa>0$ be such that $L\geq 2$ and $\tau < 1/(10L)$. Then, there exists a constant $b_0 = b_0(L,\tau,\kappa)$ satisfying the following property: If $u : \mathbb{Z}^2 \to \mathbb{R}$ is a non-constant harmonic function with $\overline{\delta}^-(u) < \sqrt{2}-1$ for which there exist $A\geq 1$ and $\widetilde{R_0} \in \mathbb{Z}$ satisfying
    \begin{equation}\label{eq:ThreeSquaresIneq}
        M_u(Q_{LR}) \leq A e^{\tau d} M_u(Q_R) + A e^{-(2-\tau)d} M_u(Q_{L^2 R})
    \end{equation}
    for any $R\geq \widetilde{R}_0$ and $d \in [1,\kappa R] \cap \mathbb{Z}$, then there is $\widetilde{R}_1 = \widetilde{R}_1(u)$ such that
    \begin{equation*}
        \|u\|_{Q_R} \geq e^{b_0 R}
    \end{equation*}
    for all $R\geq \widetilde{R}_1$.
\end{lemma}

\begin{proof}
    For each $j \in \mathbb{N}$, let 
    \begin{equation*}
        R_j = L^j \quad \text{and} \quad h_j = \log(M_u(Q_{R_{j+1}})) - \log(M_u(Q_{R_j})) = \log\left( \frac{M_u(Q_{L R_j})}{M_u(Q_{R_j})} \right).
    \end{equation*}
    Since $Q_{R_j} \subseteq Q_{R_{j+1}}$, it follows that $M_u(Q_{R_{j+1}}) \geq M_u(Q_{R_{j}})$. As a result, $(h_j)$ is a sequence of non-negative numbers. 

    To begin with, we shall prove that $(h_j)$ is unbounded. For the sake of contradiction, let us assume that $(h_j)$ is bounded. Hence, there exists $K$ such that $h_j\leq K$ for each $j \in \mathbb{N}$. Since $R_0 = 1$, we get
    \begin{equation*}
        \log(M_u(Q_{L^n})) = \log M_u(Q_1) + \sum_{j=0}^{n-1} h_j \leq \log M_u(Q_1) + Kn
    \end{equation*}
    for each $n \in \mathbb{N}$. Let $R\geq 1$. Then, there exists $n\in \mathbb{N}$ such that $L^n\leq R < L^{n+1}$. It follows that $n\leq \log(R)/\log(L)$. As a result,
    \begin{equation*}
        M_u(Q_R) \leq M_u(Q_{L^{n+1}}) \leq M_u(Q_1) e^{K(n+1)} \leq M_u(Q_1) e^K R^{K/\log(L)} = C_u R^{C_u'}.
    \end{equation*}
    This means that $u$ has polynomial growth, which contradicts Theorem \ref{thm:polynomial-main} because $\overline{\delta}^-(u) < \sqrt{2}-1$. Therefore, $(h_j)$ is unbounded.

    For each $j \in \mathbb{N}$, let us define
    \begin{equation*}
        s_j = \min\left\{ \frac{h_j}{2\tau} , \kappa R_j \right\} .
    \end{equation*}
    Moreover, let
    \begin{equation*}
        H = \max\left\{2\log(2A) , 2 + \frac{2\log(2A)}{2-\tau}\right\} \quad \text{ and } \quad R_{\ast} = \max\left\{\widetilde{R}_0 , \frac{1}{\kappa} \left( 2 + \frac{2\log(2A)}{2-\tau} \right)\right\}. 
    \end{equation*}
    Since $(h_j)$ is unbounded and $(R_j)$ is strictly increasing, there exists $j_0 \in \mathbb{Z}^+$ such that $h_{j_0}\geq H$ and $R_{j_0}\geq \max\{R_{\ast},2H/[(2-\tau)\kappa]\}$. Then \eqref{eq:ThreeSquaresIneq} implies
    \begin{equation*}
        M_u(Q_{R_{j_0+1}}) \leq A e^{\tau d} M_u(Q_{R_{j_0}}) + A e^{-(2-\tau)d} M_u(Q_{R_{j_0+2}}).
    \end{equation*}
    Dividing both sides by $M_u(Q_{R_{j_0+1}})$, we deduce
    \begin{equation}\label{eq:LargeSquares}
        1 \leq A e^{\tau d - h_{j_0}} + A e^{h_{j_0+1} - (2 - \tau)d}.
    \end{equation}
    Since $d\leq s_{j_0}\leq h_{j_0}/(2\tau)$ and $h_{j_0} \geq 2\log(2A)$, it follows that 
    \begin{equation}\label{eq:LargeSquaresTerm1}
        A e^{\tau d - h_{j_0}} \leq A e^{-h_{j_0}/2} \leq \frac{1}{2}.
    \end{equation}
    Combining \eqref{eq:LargeSquares} and \eqref{eq:LargeSquaresTerm1}, we obtain
    \begin{equation*}
        A e^{h_{j_0+1} - (2 - \tau)d} \geq \frac{1}{2}.
    \end{equation*}
    Thus,
    \begin{equation*}
        h_{j_0+1} \geq (2 - \tau)d - \log(2A) \geq (2-\tau)(s_{j_0}-1) - \log(2A).
    \end{equation*}
    Since $s_{j_0} \geq 2 + 2\log(2A)/(2-\tau)$, we get
    \begin{equation*}
        \frac{(2-\tau) s_{j_0}}{2} \geq 2 - \tau + \log(2A),
    \end{equation*}
    which implies
    \begin{equation*}
        (2-\tau)(s_{j_0}-1) - \log(2A) \geq  \frac{(2-\tau) s_{j_0}}{2}.
    \end{equation*}
    As a result,
    \begin{equation*}
        h_{j_0+1} \geq \frac{(2-\tau) s_{j_0}}{2} \geq \min\left\{\frac{2-\tau}{4\tau} h_{j_0} , \frac{(2-\tau)\kappa}{2} R_{j_0}\right\} \geq H.
    \end{equation*}
    Since $R_{j_0+1}\geq R_{j_0}$, we can iterate the previous argument to show that
    \begin{equation*}
        h_{j+1} \geq \min\left\{\frac{2-\tau}{4\tau} h_{j} , \frac{(2-\tau)\kappa}{2} R_{j}\right\} \geq H.
    \end{equation*}
    for all $j\geq j_0$. Because $\tau < 1/(10 L)$, it follows that
    \begin{equation*}
        \frac{h_{j_0+n}}{R_{j_0+n}} \geq \min\left\{\left(\frac{2-\tau}{4\tau L}\right)^n \frac{h_{j_0}}{R_{j_0}} , \frac{(2-\tau)\kappa}{2L} \right\}
    \end{equation*}
    for each $n \in \mathbb{N}^+$. Since $h_{j_0} > 0$, there exists $J'$ large enough such that 
    \begin{equation*}
        h_j \geq \frac{(2-\tau)\kappa}{2L} R_j
    \end{equation*}
    for every $j\geq J'$. Let $R \in [R_{J'+1},\infty)\cap \mathbb{Z}$ be fixed. Then, there exists $j> J'$ such that $R_j\leq R < R_{j+1}$. Since $M_u(Q_{R_{j-1}}) \geq 1$, we deduce
    \begin{equation*}
        \log(M_u(Q_R)) \geq \log(M_u(Q_{R_j})) = \log(M_u(Q_{R_{j-1}})) + h_{j-1} \geq \frac{(2-\tau)\kappa}{2L} R_{j-1} > \frac{(2-\tau)\kappa}{2L^3} R.
    \end{equation*}
    Letting $b_0 = (2-\tau)\kappa/(4L^3)$, we conclude
    \begin{equation*}
        \|u\|_{Q_R} \geq M_u(Q_R) - 1 \geq e^{2b_0R}-1 \geq e^{b_0 R} .
    \end{equation*}
\end{proof}

Our goal now is to establish \eqref{eq:ThreeSquaresIneq}. We shall see that this estimate can be deduced from the next two lemmas. The first one gives a convenient polynomial approximation of a discrete harmonic function that allows us to estimate the maximum of $u$ in $Q_{3LR}$ by a factor with exponential decay of its maximum in $Q_{L^2 R}$. The desired approximation is obtained by means of a particular basis of the space $\mathbb{R}^{[-3LR,3LR] \cap \mathbb{Z}}$ that remains invariant under the discrete Laplacian. 

\begin{lemma}\label{lemma:Step1Lower}
    There exist $L\geq 4$ and $C_0>0$ such that if $u:\mathbb{Z}^2 \to \mathbb{R}$ is harmonic, $R \in \mathbb{Z}^+$, and $d \in [1,R] \cap \mathbb{Z}$, then there exists a discrete harmonic polynomial $P_d$ with degree at most $d$ such that
    \begin{equation}\label{eq:GrowthCondPd}
        \|u - P_d\|_{Q_{3LR}} \leq C_0 e^{-2d} M_u(Q_{L^2 R}).
    \end{equation}
\end{lemma}

\begin{proof}
    For each $v \in \{\pm e_1, \pm e_2\}$, let us consider the operator $\nabla_{v} g(x) = g(x + v) - g(x)$. It is well-known that if $g$ is a discrete harmonic function in $Q_m$ and $q,q' \in Q_{m/3}$ are neighbors, then
    \begin{equation}\label{eq:NeighborsEstimate}
        |g(q) - g(q')| \leq \frac{C}{m}\|g\|_{Q_m}
    \end{equation}
    (see \cite[eq. 15]{BuhovskyLogunovMalinnikovaSodin2022} or \cite[Theorem 14]{Heilbronn1949}). Let $S \in \mathbb{Z}^+$, $k \in [1,S/4] \cap \mathbb{Z}$, and $v_1,\dots,v_k \in \{\pm e_1, \pm e_2\}$. If $m = \lfloor S/k\rfloor$, then $m\geq 4$ and $km \leq S$. Notice that if $R \in \mathbb{N}$ and $x \in Q_R$, then $x,x+v_1 \in x + Q_{m/3}$, and thus \eqref{eq:NeighborsEstimate} yields
    \begin{equation*}
        |\nabla_{v_1} g(x)| \leq \frac{C_0}{m}\|g\|_{x + Q_m} \leq \frac{C_0}{m} \|g\|_{Q_{R+m}}.
    \end{equation*}
    Thus,
    \begin{equation}\label{eq:GradientEstimate}
        \|\nabla_{v_1} g\|_{Q_R} \leq \frac{C_0}{m} \|g\|_{Q_{R+m}}.
    \end{equation}
    Since $u$ is harmonic and the discrete Laplacian commutes with translations, letting $R=0$ in the above inequality, we get
    \begin{equation*}
        |(\nabla_{v_1} \dots \nabla_{v_k} u)(0)| \leq \frac{C_0}{m} \|\nabla_{v_2} \dots \nabla_{v_k} u\|_{Q_{m}}. 
    \end{equation*}
    Using \eqref{eq:GradientEstimate} with $R=m$, we have
    \begin{equation*}
        |(\nabla_{v_1} \dots \nabla_{v_k} u)(0)| \leq \left(\frac{C_0}{m}\right)^2 \|\nabla_{v_3} \dots \nabla_{v_k} u\|_{Q_{2m}}.
    \end{equation*}
    Iterating \eqref{eq:GradientEstimate} a further $k-2$ times and using the facts that $Q_{km}\subseteq Q_S$ and $m\geq S/(2k)$, we obtain
    \begin{equation}\label{eq:EstimateFirstDiff}
        |(\nabla_{v_1} \dots \nabla_{v_k} u)(0)| \leq \left(\frac{C_0}{m}\right)^k \|u\|_{Q_{S}} \leq \left( \frac{2C_0k}{S} \right)^k\|u\|_{Q_{S}}
    \end{equation}
    for any $S \in \mathbb{Z}^+$,  $k \in [1,S/4] \cap \mathbb{Z}$, and $v_1,\dots,v_k \in \{\pm e_1, \pm e_2\}$.

    Next, for each $k \in \mathbb{N}$, let 
    \begin{equation*}
        q_k(x) = \binom{x + \lfloor k/2 \rfloor}{k}.
    \end{equation*}
    Let $Tf(x) = f(x+1)$, $\delta_+ = T - I$, $\delta_- = I - T^{-1}$, and $\mathcal{L} = T + T^{-1} - 2I = \delta_+ \delta_-$. Notice that $\mathcal{L}$ is the one-dimensional discrete Laplacian. Moreover, if $k\geq 2$, then Pascal's identity yields
    \begin{equation*}
        \delta_- \binom{x + \lfloor k/2 \rfloor}{k} = \binom{x + \lfloor k/2 \rfloor}{k} - \binom{x-1 + \lfloor k/2 \rfloor}{k} = \binom{x-1+\lfloor k/2 \rfloor}{k-1}
    \end{equation*}
    and
    \begin{equation*}
        \delta_+ \binom{x-1+\lfloor k/2 \rfloor}{k-1} = \binom{x+\lfloor k/2 \rfloor}{k-1} - \binom{x-1+\lfloor k/2 \rfloor}{k-1} = \binom{x-1+\lfloor k/2 \rfloor}{k-2}.
    \end{equation*}
    Observe that $\lfloor (k-2)/2 \rfloor = \lfloor k/2 \rfloor - 1$, which implies
    \begin{equation*}
        \mathcal{L}q_k(x) = \delta_+ \delta_- \binom{x + \lfloor k/2 \rfloor}{k} = \binom{x+\lfloor (k-2)/2 \rfloor}{k-2} = q_{k-2}(x).
    \end{equation*}
    By letting $q_k=0$ for each $k \in \mathbb{Z}^-$, we can extend the previous identity to every $k \in \mathbb{Z}$ (since linear and constant polynomials vanish after applying $\mathcal{L}$). For each $j \in \mathbb{N}$, let
    \begin{equation*}
        D_{2j} = \mathcal{L}^j \quad \text{and} \quad D_{2j+1} = (I-T^{-1})\mathcal{L}^j.
    \end{equation*}
    Notice that if $i = 2j$ for some $j\in\mathbb{N}$, then $D_i q_k(0) = q_{k-2j}(0)$. If $k < 2j$, then $q_{k-2j} \equiv 0$. If $k>2j$, then $x$ shows up in the numerator of the fraction defining $q_{k-2j}$, which implies that $q_{k-2j}(0) = 0$. On the other hand, if $i = 2j+1$ for some $j \in \mathbb{N}$, then $D_i q_k(0) = \delta_- q_{k-2j}(0) = q_{k-2j}(0) - q_{k-2j}(-1)$. If $k > 2j+1$, then both $x$ and $x+1$ show up in the numerator of the fraction defining $q_{k-2j}$, which implies that $q_{k-2j}(0) = q_{k-2j}(-1) = 0$. If $k < 2j+1$, then $q_{k-2j}$ is constant, and thus $q_{k-2j}(0) - q_{k-2j}(-1) = 0$. If $k = 2j+1$, then $q_{k-2j}(x) = x$, which means that $q_{k-2j}(0) - q_{k-2j}(-1) = 1$. Therefore,
    \begin{equation}\label{eq:ActionDiBasis}
        D_i q_k(0) = \begin{cases}
            1, &\text{if } i=k,\\
            0, &\text{if } i\neq k.
        \end{cases}
    \end{equation}

    Now we shall show that $\{q_k\}_{k=0}^{2N}$ is a basis for the space $V = \mathbb{R}^{I_N}$ of real-valued functions defined on $I_N = [-N,N] \cap \mathbb{Z}$.

    Notice that $\{\mathds{1}_{\{-N\}},\dots,\mathds{1}_{\{N\}}\}$ spans $V$, which means that $\dim V\leq 2N+1$. On the other hand, let $\alpha_0,\dots,\alpha_{2N} \in \mathbb{R}$ be such that $\sum_{k=0}^{2N} \alpha_k q_k (x) = 0$ for all $x \in I_N$. If $i \in \{0,\dots,2N\}$, we can take $D_i$ and then set $x=0$ to get that $\alpha_i = 0$. This implies that $\{q_0,\dots,q_{2N}\}$ is linearly independent in $V$. Therefore, $\dim V = 2N+1$ and $\{q_0,\dots,q_{2N}\}$ is a basis of $V$. By \eqref{eq:ActionDiBasis}, we get that
    \begin{equation*}
        f(x) = \sum_{k=0}^{2N} (D_kf)(0) q_k(x)
    \end{equation*}
    for every $f \in V$. As a result,
    \begin{equation*}
        u(x,y) = \sum_{i,j=0}^{2N}a_{i,j} q_i(x) q_j(y)
    \end{equation*}
    for every $(x,y) \in Q_N$, where $a_{i,j} = (D_i^x D_j^y u)(0,0)$.

    Next, let us define
    \begin{equation*}
        P_d(x,y) = \sum_{i+j\leq d}a_{i,j} q_i(x) q_j(y).
    \end{equation*}
    Since $q_iq_j$ has degree $i+j$, it follows that $\deg P_d \leq d$. It remains to check that $P_d$ is harmonic and satisfies \eqref{eq:GrowthCondPd}.

    For harmonicity, observe that $D_{i+2} = D_i \mathcal{L}$ and $(\mathcal{L}_x + \mathcal{L}_y)u \equiv 0$. Thus,
    \begin{equation*}
        \begin{aligned}
            (\mathcal{L}_x + \mathcal{L}_y) P_d &= \sum_{i+j\leq d}a_{i,j}(q_{i-2}(x) q_j(y) + q_i(x)q_{j-2}(y))\\ &= \sum_{i+j\leq d-2} (a_{i+2,j} + a_{i,j+2})q_i(x)q_j(y)\\ &= \sum_{i+j\leq d-2} (D_i^x D_j^y(\mathcal{L}_x + \mathcal{L}_y)u)(0,0) q_i(x)q_j(y)\\ &\equiv 0.
        \end{aligned}
    \end{equation*}

    Next, let $R\geq 1$ and $d \in [1,R]\cap \mathbb{Z}$ be fixed. Let us consider $N = 3LR$. Hence, 
    \begin{equation}\label{eq:DecomU-Pd}
        u-P_d = \sum_{\substack{0\leq i,j\leq 2N\\ i+j>d}} a_{i,j}q_i(x) q_j(y)
    \end{equation}
    in $Q_N$. Let $L\geq 48$. It follows that $4N = 12LR \leq L^2R/4$. If $(i,j)$ is one of the pairs in the sum \eqref{eq:DecomU-Pd} and $k = i+j$, then $1\leq k\leq 4N$ and $D_i^x D_j^y$ is the composition of exactly $k$ first-difference operators. By \eqref{eq:EstimateFirstDiff}, we get
    \begin{equation*}
        |a_{i,j}| \leq \left( \frac{Ck}{L^2R} \right)^k \|u\|_{Q_{L^2 R}}.
    \end{equation*}
    Notice that if $|x|\leq N$ and $i \in \{1,\dots,2N\}$, then
    \begin{equation*}
        \begin{aligned}
            |q_i(x)| = \frac{1}{i!} \prod_{\ell=0}^{i-1}|x+\lfloor i/2\rfloor - \ell| \leq \frac{(2N)^i}{i!}.
        \end{aligned}
    \end{equation*}
    This estimate also holds when $i=0$. Therefore,
    \begin{equation*}
    \begin{aligned}
        \left| \sum_{\substack{0\leq i,j\leq 2N\\ i+j=k}} a_{i,j} q_i(x) q_j(y) \right| &\leq \left( \frac{2NCk}{L^2R} \right)^k \left( \sum_{i+j=k} \frac{1}{i!j!} \right) \|u\|_{Q_{L^2 R}}\\
        &= \frac{2^k}{k!}\left( \frac{2NCk}{L^2R} \right)^k\|u\|_{Q_{L^2 R}}\\ &\leq \left( \frac{4eNC}{L^2R} \right)^k\|u\|_{Q_{L^2 R}}\\ &=  \left( \frac{12eC}{L} \right)^k\|u\|_{Q_{L^2 R}}.
    \end{aligned}
    \end{equation*}
    for all $(x,y) \in Q_N$. Letting $L\geq 48$ be large enough such that $12eC/L \leq e^{-2}$, we conclude
    \begin{equation*}
    \begin{aligned}
        \|u - P_d\|_{Q_{3 L R}} &\leq  \left( \sum_{k=d+1}^{4N} e^{-2k} \right)\|u\|_{Q_{L^2 R}}\\ &\leq \frac{e^{-2d}}{e^2-1}\|u\|_{Q_{L^2 R}}.
        \end{aligned}
    \end{equation*}
\end{proof}

Our second main ingredient to prove \eqref{eq:ThreeSquaresIneq} is the following lemma, which states that if a discrete harmonic polynomial takes negative values in $Q_{3LR}$ in a set of small density, then we can estimate its maximum in $Q_{LR}$ by an exponential growth factor of its maximum in $Q_R$. 

\begin{lemma}\label{lemma:Step3Lower}
    Let $L$ be that of Lemma \ref{lemma:Step1Lower}, and let $R \in \mathbb{Z}^+$. For every $\tau > 0$ there exist $\kappa,\varepsilon,C_{\tau} > 0$ such that if $P$ is a discrete harmonic polynomial of degree at most $d \in [1,\kappa R] \cap \mathbb{Z}$ satisfying
    \begin{equation}\label{eq:HypStep3Lower}
        |\{x \in Q_{3LR} : P(x) < 0\}| \leq \varepsilon |Q_{3LR}|,
    \end{equation}
    then 
    \begin{equation*}
        \|P\|_{Q_{LR}} \leq C_\tau e^{\tau d}\|P\|_{Q_{R}}.
    \end{equation*}
\end{lemma}

The proof of Lemma \ref{lemma:Step3Lower} requires two technical results on the growth of discrete harmonic polynomials, and thus we will first see how Lemmas \ref{lemma:Step1Lower} and \ref{lemma:Step3Lower} are combined to deduce the exponential lower bound of $u$.

\begin{proof}[Proof of Theorem \ref{thm:LowerExponentialBound}]
    Let $L\geq 4$ and $C_0>0$ be the constants from Lemma \ref{lemma:Step1Lower}, and let $\tau \in (0,1/(10L))$. Then, there exist $\kappa,\varepsilon,C_{\tau} > 0$ satisfying Lemma \ref{lemma:Step3Lower}. Without loss of generality, we can assume that $\kappa\leq 1$ and $C_{\tau} \geq 1$. Let 
    \begin{equation*}
        \delta_L = \frac{1}{2}\min\{\varepsilon , \sqrt{2} - 1\},
    \end{equation*}
    and let $u : \mathbb{Z}^2 \to \mathbb{R}$ be a non-constant harmonic function with $\overline{\delta}^-(u) < \delta$. By the definition of $\overline{\delta}^-(u)$ and the choice of $\delta$, we know that there exists $S_0 = S_0(u)$ such that
    \begin{equation*}
        |Q_{R} \cap \{u < 0\}| \leq \varepsilon |Q_R|
    \end{equation*}
    for every $R\geq S_0$. Let $R_0 \in \mathbb{Z}^+$ be such that $3 L R_0 \geq S_0$, and let $R\geq R_0$ and $d \in [1,\kappa R] \cap \mathbb{Z}$ be fixed. By Lemma \ref{lemma:Step1Lower}, there exists a discrete harmonic polynomial $P_d$ of degree at most $d$ satisfying
    \begin{equation*}
        \|u - P_d\|_{Q_{3LR}} \leq C_0 e^{-2d} M_u(Q_{L^2 R}).
    \end{equation*}
    Let $B = C_0 e^{-2d} M_u(Q_{L^2 R})$ and $P = P_d + B$. Notice that 
    \begin{equation}\label{eq:InftyNormP-U}
        0\leq P(x,y) - u(x,y)\leq 2B
    \end{equation}
    for every $(x,y) \in Q_{3LR}$, which implies that $Q_{3LR} \cap \{P<0\} \subseteq Q_{3LR} \cap \{u<0\}$. Thus,
    \begin{equation*}
        |Q_{3LR} \cap \{P < 0\}| \leq |Q_{3LR} \cap \{u < 0\}| \leq \varepsilon |Q_{3LR}|.
    \end{equation*}
    By Lemma \ref{lemma:Step3Lower}, we get that 
    \begin{equation*}
        \|P\|_{Q_{LR}} \leq C_\tau e^{\tau d}\|P\|_{Q_{R}}.
    \end{equation*}
    Observe that the triangle inequality together with \eqref{eq:InftyNormP-U} imply
    \begin{equation*}
        \|P\|_{Q_R} \leq \|u\|_{Q_R} + 2B
    \end{equation*}
    and
    \begin{equation*}
        \|u\|_{Q_{LR}} \leq \|P\|_{Q_{LR}} + 2B.
    \end{equation*}
    Combining the last three inequalities, we deduce
    \begin{equation*}
    \begin{aligned}
        \|u\|_{Q_{LR}} &\leq C_\tau e^{\tau d}\left( \|u\|_{Q_R} + 2B \right) + 2B\\ &= C_\tau e^{\tau d}\|u\|_{Q_R} + 2
        B( C_\tau e^{\tau d} + 1)\\ &\leq C_\tau e^{\tau d} M_u(Q_R) + 2C_0 (C_\tau + 1) M_u(Q_{L^2 R}) e^{-(2-\tau)d}.
    \end{aligned}
    \end{equation*}
    As a result,
    \begin{equation*}
        M_u(Q_{LR}) \leq C_\tau e^{\tau d} M_u(Q_R) + 2C_0 (C_\tau + 1) M_u(Q_{L^2 R}) e^{-(2-\tau)d}.
    \end{equation*}
    Taking $A = \max\{C_{\tau} , 2C_0(C_\tau+1)\}$ and applying Lemma \ref{lemma:ReductionLower} completes the proof.
\end{proof}

We now aim to deduce Lemma \ref{lemma:Step3Lower}. First, we introduce a result that allows us to estimate the Fourier coefficients of a discrete harmonic function by a term that depends exponentially on the degree of the polynomial. 

\begin{lemma}\label{lemma:Step3.2Lower}
    Let $L\in\mathbb{Z}^+$ and let $P : \mathbb{Z}^2 \to \mathbb{R}$ be a harmonic polynomial of degree at most $d$, where $d \in [1,\kappa_0 R]$ for a sufficiently small absolute constant $\kappa_0$. If 
    \begin{equation*}
        p_{\rho}(\theta) = P(\rho \cos\theta , \rho \sin\theta),
    \end{equation*}
    then 
    \begin{equation*}
        |\widehat{p_{\rho}}(j)| \leq \left( \|P\|_{Q_R} \right) \frac{(C_L d)^{|j|}}{|j|!} \exp\left( C_L \frac{d^2}{R} \right)
    \end{equation*}
    for every $\rho \in [0,3LR]$ and $|j|\leq d$.
\end{lemma}

\begin{proof}
    Since $P$ is real-valued, it follows that $\widehat{p_{\rho}}(-j) = \overline{\widehat{p_{\rho}}(j)}$. Therefore, it suffices to prove the desired estimate for $j \in \{0,\dots,d\}$. 
    Let us consider $z = x+iy$, $\overline{z} = x-iy$,  $\partial_z = (\partial_x - i\partial_y)/2$ and $\partial_{\overline{z}} = (\partial_x + i\partial_y)/2$. It is well-known that $P$ can be expressed in the $z$ and $\overline{z}$ variables as
    \begin{equation*}
        P(x,y) = \sum_{a+b\leq d}\frac{(\partial_z^a \partial_{\overline{z}}^bP)(0,0)}{a!b!} z^a \overline{z}^b.
    \end{equation*}
    If $z = \rho e^{i\theta}$, then $z^a \overline{z}^b = \rho^{a+b} e^{i\theta(a-b)}$, which implies that
    \begin{equation}\label{eq:ExpansionFourierCoeff}
        \begin{aligned}
            \widehat{p_{\rho}}(j) &= \frac{1}{2\pi}\int_0^{2\pi} p_{\rho}(\theta)e^{-ij\theta} \, d\theta\\ &= \frac{1}{2\pi}\sum_{a+b\leq d}\frac{(\partial_z^a \partial_{\overline{z}}^bP)(0,0)}{a!b!} \rho^{a+b}\int_0^{2\pi} e^{i\theta (a-b-j)} d\theta\\ &= \sum_{\substack{a+b\leq d\\ a-b=j}} \frac{\rho^{a+b}}{a!b!}(\partial_z^a \partial_{\overline{z}}^bP)(0,0)\\ &= \sum_{\ell=0}^{\lfloor (d-j)/2 \rfloor} \frac{\rho^{j+2\ell}}{\ell!(j+\ell)!} (\partial_z^{j+\ell}\partial_{\overline{z}}^\ell P)(0,0).
        \end{aligned}
    \end{equation}
    Thus, it remains to estimate the coefficients $(\partial_z^{j+\ell}\partial_{\overline{z}}^\ell P)(0,0)$.
    
    Let us consider the operators $\delta_x F(x,y) = F(x+1,y) - F(x,y)$ and $\delta_y F(x,y) = F(x,y+1) - F(x,y)$. Now let $t = C_1d/R$. Without loss of generality, we can choose $\kappa_0$ such that $\kappa_0\leq 1/(4C_1)$. Since $d\leq \kappa_0 R$, it follows that $d\leq R/4$ and $t\leq 1/4$. As a result, 
    \eqref{eq:EstimateFirstDiff} yields
    \begin{equation}\label{eq:BoundNormt}
        |\delta_x^a \delta_y^b P(0,0)| \leq t^{a+b} \|P\|_{Q_R}
    \end{equation}
    for all $a,b \in \mathbb{N}$ such that $a+b \leq d$. 

    For each polynomial $F$ of degree at most $d$, we define
    \begin{equation*}
        \|F\|_t = \max_{\substack{a,b \in \mathbb{N}\\ a+b\leq d}} t^{-(a+b)} |\delta_x^a \delta_y^b F(0,0)|.
    \end{equation*}
    Notice that if $a+b < d$, then
    \begin{equation*}
        t^{-(a+b)} |\delta_x^{a+1}\delta_y^b F(0,0)| \leq t \|F\|_t.
    \end{equation*}
    Moreover, if $a+b=d$, then the left-hand side of the above inequality vanishes because $\delta_x^{a+1}\delta_y^b$ has order $d+1$. As a result, $\|\delta_x F\|_t\leq t \|F\|_t$. Similarly, we deduce that $\|\delta_y F\|_t\leq t \|F\|_t$.
    Next, we shall prove that
    \begin{equation}\label{eq:DecompositionFTrans}
        F(x+s,y) = \sum_{m=0}^d \binom{s}{m} \delta_x^m F(x,y)
    \end{equation}
    for every $s \in \mathbb{R}$. Observe that if $s = n \in \mathbb{N}$, then $F(x+n,y) = (I + \delta_x)^nF(x,y)$, and thus the binomial theorem implies \eqref{eq:DecompositionFTrans} for every $s \in \mathbb{N}$. Since both sides of \eqref{eq:DecompositionFTrans} are polynomials in the $s$-variable and they coincide at infinitely many points, we deduce that \eqref{eq:DecompositionFTrans} holds for every $s \in \mathbb{R}$. 

    Observe that
    \begin{equation*}
        \frac{d}{ds}\binom{s}{m}\Big|_{s=0} = \frac{(-1)(-2)\dots(-(m-1))}{m!} = \frac{(-1)^{m-1}}{m}
    \end{equation*}
    for all $m\geq 1$. Hence, by differentiating both sides with respect to $s$ of \eqref{eq:DecompositionFTrans} and letting $s=0$, we get
    \begin{equation*}
        \partial_x F = \sum_{m=1}^d \frac{(-1)^{m-1}}{m}\delta_x^m F.
    \end{equation*}
    Applying the triangle inequality for $\|\cdot\|_t$, iterating $\|\delta_x F\|_t\leq t \|F\|_t$, and using the fact that $t\leq 1/4$, we have
    \begin{equation*}
        \|\partial_x F\|_t \leq \sum_{m=1}^d \frac{t^m}{m}\|F\|_t \leq \frac{t}{1-t}\|F\|_t \leq 2t \|F\|_t.
    \end{equation*}
    Moreover, we can perform a similar argument in the $y$-variable to deduce that $\|\partial_y F\|_t \leq 2t \|F\|_t$. This implies that $\|\partial_z F\|_t\leq 2t \|F\|_t$ and $\|\partial_{\overline{z}} F\|_t\leq 2t \|F\|_t$.

    Let us suppose now that $F$ is also discrete harmonic. This means that
    \begin{equation}\label{eq:ExpandDiscreteHarmonic}
        F(x,y) = \frac{1}{4} [F(x+1,y) + F(x-1,y) + F(x,y+1) + F(x,y-1)]
    \end{equation}
    for each $(x,y) \in \mathbb{Z}^2$. Since all the terms in the above equation are polynomials in the $x$ and $y$ variables, it follows that \eqref{eq:ExpandDiscreteHarmonic} holds for every $(x,y) \in \mathbb{R}^2$. By adding up the Taylor expansions of $F(x+1,y)$ and $F(x-1,y)$, we obtain
    \begin{equation*}
        F(x+1,y) + F(x-1,y) - 2F(x,y) = \partial_x^2 F + 2\sum_{m=2}^{\lfloor d/2 \rfloor} \frac{\partial_x^{2m} F}{(2m)!}.
    \end{equation*}
    We can also get the analogous equation in the $y$-variable. Adding up the two equations and using \eqref{eq:ExpandDiscreteHarmonic}, we deduce
    \begin{equation*}
        (\partial_x^2 + \partial_y^2)F = -2\sum_{m=2}^{\lfloor d/2 \rfloor} \frac{(\partial_x^{2m} + \partial_y^{2m}) F}{(2m)!}.
    \end{equation*}
    Since $4 \partial_z \partial_{\overline{z}} = \partial_x^2 + \partial_y^2$, we get
    \begin{equation*}
        \partial_z \partial_{\overline{z}} F = -\frac{1}{2} \sum_{m=2}^{\lfloor d/2 \rfloor} \frac{(\partial_x^{2m} + \partial_y^{2m}) F}{(2m)!}.
    \end{equation*}
    Iterating $\|\partial_x F\|_t \leq 2t \|F\|_t$ and $\|\partial_y F\|_t \leq 2t \|F\|_t$ and using the fact that $t\leq 1/2$, we obtain
    \begin{equation}\label{eq:EstimateDerivativeHarmonic}
        \|\partial_z \partial_{\overline{z}} F\|_t \leq \sum_{m=2}^{\lfloor d/2 \rfloor} \frac{(2t)^{2m}}{(2m)!} \|F\|_t \leq C t^4 \|F\|_t.
    \end{equation}
    Returning to \eqref{eq:ExpansionFourierCoeff}, observe that $(\partial_z \partial_{\overline{z}})^{\ell} P$ is discrete harmonic for each $\ell \in \mathbb{N}$, because translations and derivatives commute. By \eqref{eq:EstimateDerivativeHarmonic} and \eqref{eq:BoundNormt}, we get
    \begin{equation*}
    \begin{aligned}
         |(\partial_z^{j+\ell}\partial_{\overline{z}}^\ell P)(0,0)| &\leq \|\partial_z^j (\partial_z \partial_{\overline{z}})^{\ell} P\|_t\\  &\leq (2t)^j\|(\partial_z \partial_{\overline{z}})^{\ell} P\|_t \\ &\leq (2t)^j (Ct^4)^{\ell}\|P\|_t\\ &\leq (2t)^j (Ct^4)^{\ell} \|P\|_{Q_R}.
    \end{aligned}
    \end{equation*}
    Combining these inequalities with \eqref{eq:ExpansionFourierCoeff} and $(j+\ell)!\geq (j!)(\ell!)$, we get
    \begin{equation*}
    \begin{aligned}
        |\widehat{p_{\rho}}(j)| &\leq \left( \|P\|_{Q_R} \right) \frac{(2\rho t)^j}{j!} \sum_{\ell=0}^{\infty} \frac{(C\rho^2 t^4)^{\ell}}{(\ell!)^2}\\ &\leq \left( \|P\|_{Q_R} \right) \frac{(2\rho t)^j}{j!} \exp(2\sqrt{C} \rho t^2)\\ &\leq \left( \|P\|_{Q_R} \right) \frac{(6C_1 d L)^j}{j!} \exp(6C_1^2L\sqrt{C}d^2/R).
    \end{aligned}
    \end{equation*}
    This finishes the proof.
\end{proof}

Next, we introduce a result that estimates the $L^{\infty}$-norm of a trigonometric polynomial in terms of the size of its Fourier coefficients when the polynomial takes negative values in a set of small measure.

\begin{lemma}\label{lemma:Step3.1Lower}
    Let $p$ be a real trigonometric polynomial of degree at most $d\geq 1$ such that $p<0$ on a set of measure at most $\eta \in (0,1/128]$ (with respect to the normalized measure in the torus $\mathbb{R} / 2\pi\mathbb{Z}$). If $N = 32(d+1)$ and $s_0 = \lfloor \eta N \rfloor$, then
    \begin{equation*}
        \|p\|_{\infty} \leq C (d+1)^3 \exp(C\eta \log(e/\eta) d) \max_{|j|\leq s_0}|\widehat{p}(j)|.
    \end{equation*}
\end{lemma}

\begin{proof}
    For each $\alpha \in [0,2\pi]$ and $n \in \{0,\dots,N-1\}$, let 
    \begin{equation*}
        \theta_n(\alpha) = \alpha + \frac{2\pi n}{N}
    \end{equation*}
    and let $S(\alpha)$ be the number of these points where $p$ takes negative values, i.e. $S(\alpha) = |\{\theta_n(\alpha)\}_{n=0}^{N-1} \cap E|$. Let $E = \{p<0\}$. Then,
    \begin{equation*}
            \frac{1}{2\pi}\int_0^{2\pi}S(\alpha) \, d\alpha = \frac{1}{2\pi}\int_0^{2\pi} \left( \sum_{n=0}^{N-1} \mathds{1}_E\left(  \alpha + \frac{2\pi n }{N}\right)\right)\, d\alpha = N |E| \leq N\eta.
    \end{equation*}
    As a result, there exists $\alpha' \in [0,2\pi]$ such that $S(\alpha') \leq \lfloor N\eta \rfloor = s_0$. This means that there exist $\beta_1,\dots,\beta_{S(\alpha')} \in \{\theta_n(\alpha')\}_{n=0}^{N-1} \cap E$. Let us define
    \begin{equation*}
        w(\theta) = \prod_{\nu=1}^{S(\alpha')}|e^{i\theta} - e^{i\beta_\nu}|^2.
    \end{equation*}
    If $S(\alpha')=0$, we define $w\equiv 1$. Notice that each term in this product can be written as
    \begin{equation}\label{eq:SingleWTerm}
        |e^{i\theta} - e^{i\beta_\nu}|^2 = 2 - e^{-i\beta_\nu}e^{i\theta} - e^{i\beta_\nu}e^{-i\theta},
    \end{equation}
    which implies that $\deg w\leq S(\alpha')$. Notice that if $q(\theta) = \sum_{k} a_k e^{ik\theta}$ is a trigonometric polynomial, then
    \begin{equation*}
        q_1(\theta) = q(\theta)(2 - e^{-i\beta}e^{i\theta} - e^{i\beta}e^{-i\theta}) = \sum_k 
        (2a_k - e^{-i\beta} a_{k-1} - e^{i\beta}a_{k+1})e^{ik\theta},
    \end{equation*}
    which implies that the sum of the absolute values of the coefficients of $q_1$ is at most $4 \sum |a_k|$. Letting $q\equiv 1$ and iterating this argument $S(\alpha')$ times, we get
    \begin{equation}\label{eq:BoundWHat}
        \sum_{|k| \leq S(\alpha')} |\widehat{w}(k)| \leq 4^{S(\alpha')}.
    \end{equation}
    Let $p^+(\theta) = \max\{p(\theta),0\}$. Notice that if $p(\theta_n(\alpha')) < 0$ for some $n \in \{0,\dots,N-1\}$, then it coincides with one of the $\beta_\nu$'s, which implies that $w(\theta_n(\alpha')) = 0$, and thus
    \begin{equation}\label{eq:PositiveP}
        p^+(\theta_n(\alpha')) w(\theta_n(\alpha')) = p(\theta_n(\alpha')) w(\theta_n(\alpha')).
    \end{equation}
    This equation also holds when $p(\theta_n(\alpha')) \geq 0$, since $p^+(\theta_n(\alpha')) = p(\theta_n(\alpha'))$ in such a case. Hence, \eqref{eq:PositiveP} holds for all $n \in \{0,\dots,N-1\}$. Moreover, notice that the Fourier inversion formula implies
    \begin{equation}\label{eq:ExpFourierPW}
    \begin{aligned}
        \widehat{pw}(0) &= \frac{1}{2\pi} \int_0^{2\pi}p(\theta)w(\theta) \, d\theta\\ &= \sum_{|k|\leq S(\alpha')} \widehat{w}(k) \left( \frac{1}{2\pi} \int_0^{2\pi}p(\theta) e^{ik\theta}\,d\theta \right)\\ &=  \sum_{|k|\leq S(\alpha')} \widehat{p}(-k)\widehat{w}(k).
    \end{aligned}
    \end{equation}
    It is also well-known that if $q$ is a trigonometric polynomial of degree less than $N$ and $\lambda_1,\dots,\lambda_N$ are equally spaced points in $[0,2\pi]$ modulo $2\pi$, then
    \begin{equation}\label{eq:IdentityTrigPol}
        \frac{1}{N} \sum_{k=1}^N q(\lambda_k) = \frac{1}{2\pi}\int_0^{2\pi}q(\theta) \, d\theta = \widehat{q}(0).
    \end{equation}
    Applying this result for $pw$ (which has degree at most $d + S(\alpha') < N$) and using \eqref{eq:PositiveP}-\eqref{eq:ExpFourierPW}, we deduce
    \begin{equation*}
        \frac{1}{N} \sum_{n=0}^{N-1}p^+(\theta_n(\alpha')) w(\theta_n(\alpha')) = \widehat{pw}(0) = \sum_{|k|\leq S(\alpha')} \widehat{p}(-k)\widehat{w}(k).
    \end{equation*}
    By taking absolute value on both sides, applying the triangle inequality, and using \eqref{eq:BoundWHat}, we get
    \begin{equation}\label{eq:IneqP+}
        \frac{1}{N} \sum_{n=0}^{N-1}p^+(\theta_n(\alpha')) w(\theta_n(\alpha')) \leq 4^{S(\alpha')} \max_{|j|\leq s_0}|\widehat{p}(j)|.
    \end{equation}
    Let us fix $n \in \{0,\dots,N-1\}$ such that $p(\theta_n(\alpha'))\geq 0$. For each $\nu \in \{1,\dots,S(\alpha')\}$, let $r_{\nu}$ be the number of grid steps from $\beta_\nu$ to $\theta_n(\alpha')$ in the grid $\{\theta_n(\alpha')\}_{n=0}^{N-1}$, measured along the shorter arc of the circle. We can reorder the $\beta_\nu$'s to get
    \begin{equation*}
        1\leq r_1\leq r_2\leq \dots \leq r_{S(\alpha')} \leq N/2.
    \end{equation*}
    Notice that there are at most two grid points at each positive distance from $\theta_n(\alpha')$. Thus, there are at most $2r_\nu$ grid points at distances $1,\dots,r_\nu$. Since all the first $\nu$ points in $\{\theta_n(\alpha')\}_{n=0}^{N-1} \cap E$ lie among them, we have that $\nu \leq 2r_\nu$. Moreover, notice that the angle between $\theta_n(\alpha')$ and $\beta_\nu$ (viewed as points on the unit circle) is $2\pi r_\nu/N$. By \eqref{eq:SingleWTerm} and the fact that $\sin(x)\geq 2x/\pi$ for each $x \in [0,\pi/2]$, we get
    \begin{equation*}
        |e^{i\theta_n(\alpha')} - e^{i\beta_\nu}|^2 = 2 - 2\cos\left( \frac{2\pi r_\nu}{N} \right) = 4 \sin^2 \left( \frac{\pi r_\nu}{N} \right) \geq \frac{4 \nu^2}{N^2}.
    \end{equation*}
    Thus,
    \begin{equation*}
        w(\theta_n(\alpha')) \geq \prod_{\nu=1}^{S(\alpha')}\frac{4 \nu^2}{N^2} = \left( \frac{4}{N^2} \right)^{S(\alpha')} (S(\alpha')!)^2.
    \end{equation*}
    Since $p^+(\theta) = 0$ at every point in $\{\theta_n(\alpha')\}_{n=0}^{N-1} \cap E$, then \eqref{eq:IneqP+} and the above inequality yield
    \begin{equation*}
        \frac{1}{N} \sum_{n=0}^{N-1}p^+(\theta_n(\alpha')) \leq \left( \frac{N^{S(\alpha')}}{S(\alpha')!} \right)^2 \max_{|j|\leq s_0}|\widehat{p}(j)|.
    \end{equation*}
    Notice that the inequality $S(\alpha')^{S(\alpha')}/S(\alpha')! \leq \exp(S(\alpha'))$, the monotonicity of $x\mapsto x \log(e/x)$ in $(0,1]$, and the fact that $S(\alpha')/N \leq \eta$ yield
    \begin{equation*}
        \log\frac{N^{S(\alpha')}}{S(\alpha')!} \leq S(\alpha') \log \frac{eN}{S(\alpha')} \leq N\eta \log(e/\eta).
    \end{equation*}
    Thus,
    \begin{equation*}
        \frac{1}{N} \sum_{n=0}^{N-1}p^+(\theta_n(\alpha')) \leq \exp(2N\eta \log(e/\eta))\max_{|j|\leq s_0}|\widehat{p}(j)|.
    \end{equation*}
    Since $p$ is real-valued, we know that $|p| = 2p^+ - p$. By using the above inequality and \eqref{eq:IdentityTrigPol} with $q=p$, we get
    \begin{equation}\label{eq:IneqAbsP}
    \begin{aligned}
        \frac{1}{N} \sum_{n=0}^{N-1}|p(\theta_n(\alpha'))| &\leq [2\exp(2N\eta \log(e/\eta))+1] \max_{|j|\leq s_0}|\widehat{p}(j)|\\ &\leq 3\exp(2N\eta \log(e/\eta))\max_{|j|\leq s_0}|\widehat{p}(j)|.
    \end{aligned}
    \end{equation}
    Observe that for any $|k|\leq d$, the trigonometric polynomial $q(\theta) = p(\theta) e^{-ik\theta}$ has degree at most $2d < N$, and thus an application of \eqref{eq:IdentityTrigPol} yields
    \begin{equation}\label{eq:FourierCoeffEquid}
        \widehat{p}(k) = \frac{1}{N}\sum_{n=0}^{N-1} p(\theta_n(\alpha')) e^{-ik\theta_n(\alpha')}.
    \end{equation}
    Using the Fourier inversion formula, \eqref{eq:IneqAbsP}, and \eqref{eq:FourierCoeffEquid}, we deduce
    \begin{equation*}
    \begin{aligned}
        \|p\|_{\infty} &\leq 3(2d+1)\exp(2N\eta \log(e/\eta))\max_{|j\leq s_0}|\widehat{p}(j)|\\ &\leq 3(2d+1)\exp(128\eta \log(e/\eta) d)\max_{|j\leq s_0}|\widehat{p}(j)|.
    \end{aligned}
    \end{equation*}
    This finishes the proof.
\end{proof}

Finally, we combine Lemmas \ref{lemma:Step3.2Lower} and \ref{lemma:Step3.1Lower} to obtain Lemma \ref{lemma:Step3Lower}, which finishes the proof of Theorem \ref{thm:LowerExponentialBound}.

\begin{proof}[Proof of Lemma \ref{lemma:Step3Lower}]
   Let $E = [-3LR,3LR]^2 \cap \{P<0\}$. To begin with, we shall prove that 
   \begin{equation}\label{eq:EstimateAreaElattice}
       ||E| - |E\cap \mathbb{Z}^2|| \leq 2(d+1)(12LR+1).
   \end{equation}
   Let $y \in [-3LR,3LR]$ be fixed. Observe that $x\mapsto P(x,y)$ is a polynomial in the $x$-variable of degree at most $d$. If this polynomial does not identically vanish, then it has at most $d$ roots. As a result, the set $E_y=\{x \in [-3LR,3LR] : P(x,y)<0\}$ is composed of at most $d+1$ intervals. Notice that $|E_y| = \int_{-3LR}^{3LR}\mathds{1}_E(x,y)\,dx$ and $|E_y\cap \mathbb{Z}| = \sum_{j=-3LR}^{3LR}\mathds{1}_{E}(j,y)$. Since the length of an interval and the number of integers in such an interval differ at most by $2$, it follows that
   \begin{equation}\label{eq:EstimateAreaElattice1}
       \left| \int_{-3LR}^{3LR}\mathds{1}_E(x,y)\,dx - \sum_{j=-3LR}^{3LR}\mathds{1}_{E}(j,y) \right| \leq 2(d+1).
   \end{equation}
   Notice that if  $x\mapsto P(x,y)$ is identically zero, then both terms on the left-hand side of \eqref{eq:EstimateAreaElattice1} equal zero, and thus this estimate also holds in this case. Let us consider 
   \begin{equation*}
        B = \int_{-3LR}^{3LR}\sum_{j=-3LR}^{3LR}\mathds{1}_{E}(j,y) \, dy.    
   \end{equation*}
   Since $|E| = \int_{-3LR}^{3LR}\int_{-3LR}^{3LR}\mathds{1}_E(x,y)\,dx\, dy$, we get
   \begin{equation}\label{eq:EstimateAreaElattice2}
       ||E| - B| \leq \int_{-3LR}^{3LR}\left| \int_{-3LR}^{3LR}\mathds{1}_E(x,y)\,dx - \sum_{j=-3LR}^{3LR}\mathds{1}_{E}(j,y) \right|\,dy \leq 2(d+1)(6LR).
   \end{equation}
   By changing the roles of $x$ and $y$ in the argument of \eqref{eq:EstimateAreaElattice1}, we can deduce that
   \begin{equation*}
       \left| \int_{-3LR}^{3LR}\mathds{1}_E(x,y)\,dy - \sum_{j=-3LR}^{3LR}\mathds{1}_{E}(x,j) \right| \leq 2(d+1)
   \end{equation*}
   for every $x \in [-3LR,3LR]$. Thus,
   \begin{equation}\label{eq:EstimateAreaElattice3}
       |B - |E\cap \mathbb{Z}^2|| \leq \sum_{k=-3LR}^{3LR}\left| \int_{-3LR}^{3LR}\mathds{1}_E(k,y)\,dy - \sum_{j=-3LR}^{3LR}\mathds{1}_{E}(k,j) \right| \leq 2(d+1) (6LR + 1).
   \end{equation}
   Combining the triangle inequality with \eqref{eq:EstimateAreaElattice2} and \eqref{eq:EstimateAreaElattice3}, we obtain \eqref{eq:EstimateAreaElattice}.

   Now, since $1\leq d\leq \kappa R$, we have that $(d+1)/R \leq 2 \kappa$. Moreover, observe that \eqref{eq:HypStep3Lower} implies that $|E \cap \mathbb{Z}^2| \leq \varepsilon (6LR+1)^2$. Applying \eqref{eq:EstimateAreaElattice}, we deduce
   \begin{equation}\label{eq:EstimateAreaE}
       |E| \leq \varepsilon (6LR+1)^2 + 2(d+1)(12LR+1) \leq C_L R^2 (\varepsilon + 2\kappa).
   \end{equation}
   Let $\eta \in (0,1/128]$, and let us consider $I = [2LR , 3LR]$ and $a(\rho) = |\{\theta : p_{\rho}(\theta) < 0\}|$. If $A = \{(x,y) : 2LR\leq |(x,y)| \leq 3LR\}$, then an application of polar coordinates yields
   \begin{equation*}
   \begin{aligned}
       |E| &\geq |E\cap A|\\ &= \int_{I} \rho \int_0^{2\pi} \mathds{1}_{\{p_{\rho}(\theta)<0\}} \, d\theta \, d\rho\\ &= 2\pi \int_{I} \rho \, a(\rho) \, d\rho\\ &\geq 4\pi \eta LR \cdot |I \setminus \{a(\rho) \leq \eta\}|.
   \end{aligned}
   \end{equation*}
   Dividing both sides of the previous inequality by $|I| = LR$ and using \eqref{eq:EstimateAreaE}, we get
   \begin{equation}\label{eq:Est|I-G|}
        \frac{|I \setminus \{a(\rho) \leq \eta\}|}{|I|} \leq \frac{C_L}{\eta} (\varepsilon + 2\kappa). 
   \end{equation}
   Let $\rho \in I$ be such that $a(\rho) \leq \eta$. If $a(\rho) > 0$, an application of Lemma \ref{lemma:Step3.1Lower} together with the fact that $x\mapsto x\log(e/x)$ is increasing in $(0,1)$ yield
   \begin{equation}\label{eq:boundLinftyPrho}
   \begin{aligned}
       \|p_{\rho}\|_{\infty} &\leq C(d+1)^3 \exp(C a(\rho) \log(e/a(\rho)) d) \max_{|j|\leq s_0}|\widehat{p_\rho}(j)|\\ &\leq C(d+1)^3 \exp(C \eta \log(e/\eta )d) \max_{|j|\leq s_0}|\widehat{p_\rho}(j)|
   \end{aligned}
   \end{equation}
   where $s_0 = \lfloor 32\eta (d+1) \rfloor$. If $a(\rho) = 0$, then $p_{\rho}\geq 0$, which implies that $|\widehat{p_\rho}(j)| \leq \widehat{p_{\rho}}(0)$ for every $j \in \mathbb{Z}$. By the Fourier inversion formula, it follows that $\|p_{\rho}\|_{\infty} \leq (2d+1) \widehat{p_{\rho}}(0)$. Thus, \eqref{eq:boundLinftyPrho} also holds in such a case. 

   Notice that if $j\leq s_0 \leq 64\eta d$, then the inequality $j^j/j! \leq e^j$ and the monotonicity of the function $x\mapsto x\log(eC_Ld/x)$ in $(0,C_Ld)$ yield
   \begin{equation*}
       \log \frac{(C_L d)^{j}}{j!} \leq j \log \frac{eC_Ld}{j} \leq C_L \eta \log(e/\eta)d.
   \end{equation*}
   As a consequence of \eqref{eq:boundLinftyPrho}, Lemma \ref{lemma:Step3.2Lower}, and the fact that $d^2/R\leq \kappa d$, we deduce
   \begin{equation}\label{eq:DoubleSup}
       \sup_{\{a(\rho)\leq \eta\}}\|p_{\rho}\|_{\infty} \leq C (d+1)^3 \exp(C_L[\kappa + \eta \log(e/\eta)]d)\left( \|P\|_{Q_R} \right).
   \end{equation}
   Let $\lambda \in (0,1/2]$ be such that $3\sqrt{\lambda} \leq \tau/4$, and let $\eta \in (0,1/128]$ be such that $C_L \eta \log(e/\eta) \leq \tau/4$. Moreover, let $\kappa\leq \kappa_0$ and $\varepsilon>0$ be small enough such that $C_L\kappa \leq \tau/4$ and $C_L(\varepsilon + 2\kappa)/\eta \leq \lambda$. By \eqref{eq:Est|I-G|}, we get that
   \begin{equation}\label{eq:NewDensity}
       \frac{|I \setminus \{a(\rho) \leq \eta\}|}{|I|} \leq \lambda.
   \end{equation}
   By rescaling the Remez inequality \cite[Theorem 1.1]{BrudnyiYomdin} (see also \cite[Lemma 3.1]{BuhovskyLogunovMalinnikovaSodin2022}), we have
   \begin{equation*}
       \sup_{\rho \in I}|p_{\rho}(\theta)| \leq T_d \left( \frac{2|I| - |\{a(\rho)\leq \eta\}|}{|\{a(\rho)\leq \eta\}|} \right) \sup_{\{a(\rho)\leq \eta\}}|p_{\rho}(\theta)|,
   \end{equation*}
   where $T_d$ denotes the Chebyshev polynomial of degree $d$. Observe that \eqref{eq:NewDensity} implies that $2|I| - |\{a(\rho)\leq \eta\}| \leq (1+\lambda)|I|$ and $|\{a(\rho)\leq \eta\}| \geq (1-\lambda)|I|$. Since $T_d$ is increasing in $[1,\infty)$, we deduce
   \begin{equation*}
       \sup_{\rho \in I}|p_{\rho}(\theta)| \leq T_d \left( \frac{1+\lambda}{1-\lambda} \right) \sup_{\{a(\rho)\leq \eta\}}|p_{\rho}(\theta)|.
   \end{equation*}
   Notice that there exists $t\geq 0$ such that $\cosh(t) = (1+\lambda)/(1-\lambda)$. Thus,
   \begin{equation*}
       T_d \left( \frac{1+\lambda}{1-\lambda} \right) = T_d(\cosh(t)) = \cosh(dt) \leq e^{dt}.
   \end{equation*}
   Since $t^2/2 \leq \cosh(t)-1 = 2\lambda/(1-\lambda) \leq 4\lambda$, it follows that $t\leq 3\sqrt{\lambda}$. Therefore,
   \begin{equation*}
       \sup_{\rho \in I}|p_{\rho}(\theta)| \leq e^{3\sqrt{\lambda}d}\sup_{\{a(\rho)\leq \eta\}}|p_{\rho}(\theta)|.
   \end{equation*}
   Taking the supremum with respect to $\theta$ in the above equation and using \eqref{eq:DoubleSup}, we deduce
   \begin{equation*}
       \sup_{(x,y) \in A}|P(x,y)| \leq  C (d+1)^3 \exp([3\sqrt{\lambda} +C_L\kappa + C_L\eta \log(e/\eta)]d)\left( \|P\|_{Q_R} \right).
   \end{equation*}
   Because of the choice of the constants and the fact that $C (d+1)^3 \leq C_{\tau}e^{\tau d/4}$, we obtain
   \begin{equation*}
       \sup_{(x,y) \in A}|P(x,y)| \leq C_{\tau}e^{\tau d}\left( \|P\|_{Q_R} \right).
   \end{equation*}
   Notice that if $(x,y) \in \partial Q_{2LR}$, then
   \begin{equation*}
       2LR = \|(x,y)\|_{\infty} \leq |(x,y)| \leq \sqrt{2}\|(x,y)\|_{\infty} = 2\sqrt{2}LR < 3LR.
   \end{equation*}
   This means that $\partial Q_{2LR} \subseteq A$. By the discrete maximum principle, we get
   \begin{equation*}
       \|P\|_{Q_{2LR}} \leq \sup_{(x,y) \in A}|P(x,y)| \leq C_{\tau}e^{\tau d}\left( \|P\|_{Q_R} \right).
   \end{equation*}
   This finishes the proof.
\end{proof}

\section{Exponential upper bound: Proof of Theorem \ref{thm:main}} \label{Section:ExpGrowth}
In order to establish Theorem \ref{thm:main}, we first prove Proposition \ref{prop:large-range}  by  reducing \eqref{eq:main-bound} to an inequality involving the maxima of two diamonds $D_R$. Then, we develop the argument to get the required estimates for the growth of $|u|$ on concentric translated diamonds that will lead to the final result.
\subsection{Proof of Proposition \ref{prop:large-range}}
The precise sufficient condition to obtain \eqref{eq:main-bound} is shown in the following lemma. Throughout this subsection, let \(u:\mathbb Z^2\to\mathbb R\) be
discrete harmonic, assume that \(\delta=\overline{\delta}^{-}(u)<1/32\),
and set \(\gamma=(1-32\delta)/8\).

\begin{lemma}\label{lemma:InitialReduction}
    Suppose that there exists a constant $C$, independent of $u$ and $\delta$, such that for every $R \in \mathbb{Z}^+$ large enough there exists $K \in \mathbb{Z}^+$ with
    \begin{equation*}
        \gamma R \leq K \leq \left( \frac{1}{2} - \gamma \right) R
    \end{equation*}
    satisfying
    \begin{equation}\label{eq:ReductionMainThm}
        M_u(D_{R+K}) \leq M_u(D_R) e^{C K \log(eR/K)}.
    \end{equation}
    Then, \eqref{eq:main-bound} holds for every $(x,y) \in \mathbb{Z}^2$.
\end{lemma}
\begin{proof}
   Choose $R_0$ such that the hypothesis holds for every integer $R\ge R_0$. That is, there exists $K_0 \in \mathbb{Z}^+$ with
    \begin{equation*}
        \gamma R_0 \leq K_0 \leq \left( \frac{1}{2} - \gamma \right) R_0
    \end{equation*}
    satisfying
    \begin{equation*}
        M_u(D_{R_0+K_0}) \leq M_u(D_{R_0}) e^{C K_0 \log(eR_0/K_0)}.
    \end{equation*}
    Given $R_j$, we choose $K_j$ satisfying the hypothesis at $R=R_j$, and let $R_{j+1}=R_j+K_j$. Since $\gamma R_j\leq K_j$ for each $j \in \mathbb{N}$, we deduce that 
    \begin{equation*}
        \log\left(\frac{eR_j}{K_j}\right) \leq \log\left(\frac{e}{\gamma}\right) = \log8 + \log\left( \frac{e}{1-32\delta}  \right) \leq (\log 8 + 1)\log\left( \frac{e}{1-32\delta}  \right)
    \end{equation*}
    for every $j \in \mathbb{N}$. By \eqref{eq:ReductionMainThm} and the choice of each $R_j$, we get
    \begin{equation*}
        \begin{aligned}
            \log(M_u(D_{R_j})) &\leq \log(M_u(D_{R_0})) + C (\log 8 + 1)\log\left( \frac{e}{1-32\delta}  \right) \sum_{i=0}^{j-1}K_i\\ &=\log(M_u(D_{R_0})) + C (\log 8 + 1)\log\left( \frac{e}{1-32\delta}  \right)(R_j-R_0)\\ &\leq \log(M_u(D_{R_0})) + C (\log 8 + 1)\log\left( \frac{e}{1-32\delta}  \right)R_j
        \end{aligned}
    \end{equation*}
    for all $j \in \mathbb{N}$. If $R \in \mathbb{Z}^+$ satisfies $R_j \leq R \leq R_{j+1}$ for some $j \in \mathbb{N}$, the fact that $R_{j+1} = R_j + K_j < 3R_j/2$ yields
    \begin{equation*}
        \begin{aligned}
            \log(M_u(D_R)) &\leq \log(M_u(D_{R_{j+1}}))\\ &\leq  \log(M_u(D_{R_0})) + C (\log 8 + 1)\log\left( \frac{e}{1-32\delta}  \right)R_{j+1}\\ &\leq \log(M_u(D_{R_0})) + \frac{3}{2}C (\log 8 + 1)\log\left( \frac{e}{1-32\delta}  \right)R.
        \end{aligned}
    \end{equation*}
    Therefore,
    \begin{equation*}
        M_u(D_R) \leq M_u(D_{R_0}) \exp\left( \frac{3}{2}C (\log 8 + 1)\log\left( \frac{e}{1-32\delta}  \right)R \right)
    \end{equation*}
    for every $R\geq R_0$. Letting $R = |x| + |y|$ for all $(x,y) \in \mathbb{Z}^2$ with $|x|+|y|\geq R_0$, we deduce \eqref{eq:main-bound} for all such $(x,y) \in \mathbb{Z}^2$. Because $|u(x,y)| \leq M_u(D_{R_0})$ for all the remaining $(x,y) \in \mathbb{Z}^2$, we conclude that \eqref{eq:main-bound} holds for every $(x,y) \in \mathbb{Z}^2$ with $C_u = M_u(D_{R_0})$.
\end{proof}

Now we seek to prove \eqref{eq:ReductionMainThm}. To begin with, by using the harmonicity of $u$, we shall prove that $M_u(D_{R+K})$ can be estimated by a multiple of the maximum of $|u|$ in $D_{R+K} \cap Q_{R-1}$. The constant arising in this estimate is a power of $7$, which we can join with the exponential term in \eqref{eq:ReductionMainThm} later on.

\begin{lemma}\label{lemma:Step5}
    Let $R\geq 4$ and $1\leq K\leq R/2$ be integers. Then,
    \begin{equation*}
        \|u\|_{D_{R+K}} \leq 7^{K+1} \|u\|_{D_{R+K} \cap Q_{R-1}}.
    \end{equation*}
\end{lemma}

\begin{proof}
    To begin with, we shall prove that
    \begin{equation}\label{eq:IndHypExponential}
        |u(R-1+d,y)| \leq 7^d \|u\|_{D_{R+K} \cap Q_{R-1}} 
    \end{equation}
     for any $0\leq d\leq K+1$ and $|y|\leq K+1-d$ by using induction on $d$.

    For the base case, notice that $\{(R-1,y) : |y|\leq K+1\} \subseteq D_{R+K} \cap Q_{R-1}$, which implies that
    \begin{equation*}
        |u(R-1,y)| \leq \|u\|_{D_{R+K} \cap Q_{R-1}}
    \end{equation*}
    for every $|y|\leq K+1$. This establishes the case $d=0$.

    Moreover, let $|y|\leq K$. Observe that $\{(R-2,y) : |y|\leq K+1\} \subseteq D_{R+K} \cap Q_{R-1}$. Since $u$ is harmonic, we get
    \begin{equation*}
        4u(R-1,y) = u(R,y) + u(R-2,y) + u(R-1,y+1) + u(R-1,y-1).
    \end{equation*}
    Given that $|y\pm 1| \leq K+1$, we deduce
    \begin{equation*}
        \begin{aligned}
            |u(R,y)| &\leq 4|u(R-1,y)| + |u(R-2,y)| + |u(R-1,y+1)| + |u(R-1,y-1)|\\ &\leq 7 \|u\|_{D_{R+K} \cap Q_{R-1}}.
        \end{aligned}
    \end{equation*}
    This establishes the case $d=1$.

    For the inductive step, let $d\in \{1,\dots,K\}$, and let us assume that \eqref{eq:IndHypExponential} holds for every $d'\leq d$. Let $|y| \leq K-d$. Since $|y\pm 1| \leq K+1-d$, the induction hypothesis yields that
    \begin{equation*}
        |u(R-1+d,y)|,|u(R-1+d,y\pm 1)| \leq 7^d \|u\|_{D_{R+K} \cap Q_{R-1}}
    \end{equation*}
    and
    \begin{equation*}
        |u(R-2+d,y)| \leq 7^{d-1}\|u\|_{D_{R+K} \cap Q_{R-1}}.
    \end{equation*}
    Since $u$ is harmonic, we get
    \begin{equation*}
        \begin{aligned}
            |u(R+d,y)| \leq& 4|u(R-1+d,y)| + |u(R-1+d,y+1)|\\ &+|u(R-1+d,y-1)| + |u(R-2+d,y)|\\ \leq& (6 \cdot 7^d + 7^{d-1})\|u\|_{D_{R+K} \cap Q_{R-1}} \\ \leq &7^{d+1}\|u\|_{D_{R+K} \cap Q_{R-1}}.
        \end{aligned}
    \end{equation*}
    This establishes the result for $d+1$.

    Hence, if we consider the set
    \begin{equation*}
        T_1 = \{(x,y)\in \mathbb{Z}^2 : R\leq x\leq R+K \text{ and } |y|\leq R+K-x\},
    \end{equation*}
    we have
    \begin{equation}\label{eq:BoundTriangle2}
        \|u\|_{T_1} \leq 7^{K+1}\|u\|_{D_{R+K} \cap Q_{R-1}}.
    \end{equation}
    Observe that $D_{R+K} \setminus Q_{R-1} = \bigcup_{j=1}^4 T_j$, where $T_2,T_3,T_4$ are defined by the regions illustrated in Figure \ref{fig:DecompositionDiamond}.

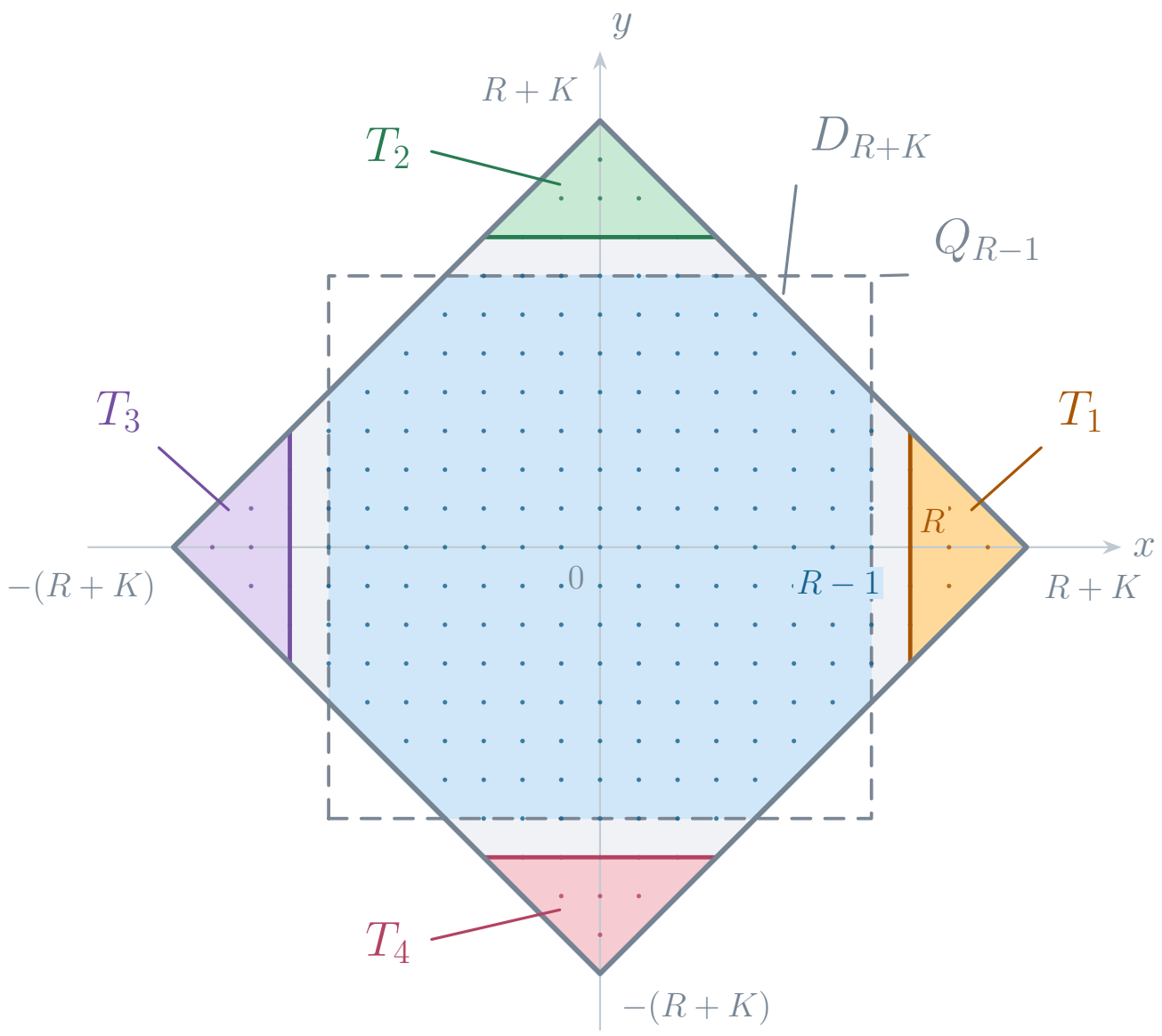
\begin{figure}[H]
    \centering

        \begin{tikzpicture}[
    x=0.48cm, y=0.48cm,
    scale=0.6, transform shape,
    line cap=round, line join=round,
    >={Stealth[length=4pt]},
    tiplabel/.style={font=\Large, inner sep=2pt},
    leader/.style={
        line width=0.6pt,
        shorten <=3pt,
        shorten >=2pt
    }
]

\def\R{8}
\def\K{3}
\pgfmathtruncatemacro{\outerR}{\R+\K}
\pgfmathtruncatemacro{\innerR}{\R-1}

\fill[GapFill]
    (\outerR,0) -- (0,\outerR) --
    (-\outerR,0) -- (0,-\outerR) -- cycle;

\begin{scope}
    \clip (-\innerR,-\innerR)
        rectangle (\innerR,\innerR);
    \fill[Efill]
        (\outerR,0) -- (0,\outerR) --
        (-\outerR,0) -- (0,-\outerR) -- cycle;
\end{scope}

\fill[TOneFill]
    (\R,\K) -- (\outerR,0) -- (\R,-\K) -- cycle;
\fill[TTwoFill]
    (-\K,\R) -- (0,\outerR) -- (\K,\R) -- cycle;
\fill[TThreeFill]
    (-\R,\K) -- (-\outerR,0) -- (-\R,-\K) -- cycle;
\fill[TFourFill]
    (-\K,-\R) -- (0,-\outerR) -- (\K,-\R) -- cycle;

\draw[
    SquareColor,
    line width=0.7pt,
    dash pattern=on 4pt off 3pt
]
    (-\innerR,-\innerR) rectangle (\innerR,\innerR);

\draw[->,AxisColor]
    ({-\outerR-2.2},0) -- ({\outerR+2.45},0)
    node[right,text=Boundary] {\large $x$};
\draw[->,AxisColor]
    (0,{-\outerR-1.45}) -- (0,{\outerR+1.8})
    node[above right,text=Boundary] {\large $y$};

\foreach \xx in {-\outerR,...,\outerR}{
    \foreach \yy in {-\outerR,...,\outerR}{
        \pgfmathtruncatemacro{\inside}{
            abs(\xx)+abs(\yy)<=\outerR
        }
        \ifnum\inside=1\relax
            \def\pointcolor{Eink}
            \ifnum\xx>\innerR\relax
                \def\pointcolor{TOneInk}
            \fi
            \ifnum\yy>\innerR\relax
                \def\pointcolor{TTwoInk}
            \fi
            \ifnum\xx<-\innerR\relax
                \def\pointcolor{TThreeInk}
            \fi
            \ifnum\yy<-\innerR\relax
                \def\pointcolor{TFourInk}
            \fi
            \fill[\pointcolor,opacity=0.8]
                (\xx,\yy) circle[radius=0.8pt];
        \fi
    }
}

\draw[TOneInk,line width=0.9pt]
    (\R,-\K) -- (\R,\K);
\draw[TTwoInk,line width=0.9pt]
    (-\K,\R) -- (\K,\R);
\draw[TThreeInk,line width=0.9pt]
    (-\R,-\K) -- (-\R,\K);
\draw[TFourInk,line width=0.9pt]
    (-\K,-\R) -- (\K,-\R);

\draw[Boundary,line width=1.1pt]
    (\outerR,0) -- (0,\outerR) --
    (-\outerR,0) -- (0,-\outerR) -- cycle;

\node[tiplabel,TOneInk,anchor=west] (T1)
    at ({\outerR+0.65},{\K+0.5}) {$T_1$};
\draw[leader,TOneInk]
    (T1.south west) -- ({\R+1.4},0.8);

\node[tiplabel,TTwoInk,anchor=east] (T2)
    at ({-\K-1.7},{\outerR-0.7}) {$T_2$};
\draw[leader,TTwoInk]
    (T2.east) -- (-0.8,{\R+1.3});

\node[tiplabel,TThreeInk,anchor=east] (T3)
    at ({-\outerR-0.65},{\K+0.5}) {$T_3$};
\draw[leader,TThreeInk]
    (T3.south east) -- ({-\R-1.4},0.8);

\node[tiplabel,TFourInk,anchor=east] (T4)
    at ({-\K-1.7},{-\outerR+0.8}) {$T_4$};
\draw[leader,TFourInk]
    (T4.east) -- (-0.8,{-\R-1.3});

\node[Boundary,font=\Large,anchor=west] (D)
    at ({\R-2.9},{\outerR-0.45}) {$D_{R+K}$};
\draw[leader,Boundary]
    (D.south west)
    -- ({\K+1.7},{\outerR-\K-1.7});

\node[SquareColor,anchor=west] (Q)
    at ({\R+0.3},{\R-0.1}) {\Large $Q_{R-1}$};
\draw[leader,SquareColor]
    (Q.south west) -- (\innerR,\innerR);

\node[below left,text=SquareColor]
    at (-0.1,-0.2) {$0$};

\node[below,Eink,fill=Efill,inner sep=1pt]
    at ({\innerR-0.85},-0.5) {$R-1$};

\node[above right,TOneInk,fill=TOneFill,inner sep=1pt]
    at ({\R+0.15},0.3) {$R$};

\node[below right,Boundary]
    at ({\outerR+0.15},-0.4) {$R+K$};

\node[below left,Boundary]
    at ({-\outerR-0.15},-0.3) {$-(R+K)$};

\node[above left,Boundary]
    at (-0.25,{\outerR+0.15}) {$R+K$};

\node[below right,Boundary]
    at (0.25,{-\outerR-0.15}) {$-(R+K)$};

\end{tikzpicture}
    \caption{Decomposition of $D_{R+K} \setminus Q_{R-1}$.}
    \label{fig:DecompositionDiamond}
\end{figure}

    Let $T : \mathbb{R}^2 \to \mathbb{R}^2$ be the counterclockwise rotation by $\pi/2$ in $\mathbb{Z}^2$. For each $j \in \{1,2,3\}$, let us consider $v_j(x,y) = u(T^j(x,y))$. It follows that each $v_j$ is a discrete harmonic function and
    \begin{equation}\label{eq:RotationsU}
        \|v_j\|_{D_{R+K} \cap Q_{R-1}} = \|u\|_{D_{R+K} \cap Q_{R-1}}
    \end{equation}
    for every $j \in \{1,2,3\}$. This implies that $\eqref{eq:IndHypExponential}$ holds for each $v_j$, which implies that \eqref{eq:BoundTriangle2} holds for $T_2,T_3,T_4$. This finishes the proof.
\end{proof}

Next, we want to estimate the maximum of $|u|$ in $D_{R+K}\cap Q_{R-1}$ by a multiple of $M_u(D_R)$. Hence, we need to handle the values of $|u|$ on the segments
\begin{equation*}
    \{(\varepsilon_1 j , \varepsilon_2(R+k-j)) : 1\leq k\leq K \text{ and } k+1 \leq j\leq R-1\}
\end{equation*}
for all $\varepsilon = (\varepsilon_1,\varepsilon_2) \in \{-1,1\}^2$. For this estimate, it will be convenient to decompose $u$ along these segments as a sum of two terms. For one of them we can establish the desired bound by only applying the harmonicity of $u$, while the other term is a polynomial, and its bound is more delicate. In fact, the assumptions on the number of negative values of $u$ are important for the polynomial estimates

\begin{lemma}\label{Lemma:4.2}
    Let $R\geq 4$ and $1\leq K\leq R/2$ be integers. For each $k \in \{1,\dots,K\}$, let us consider $I_k = \{k+1,\dots,R-1\}$. Then, there exists an absolute constant $C>0$ such that, for each $k \in \{1,\dots,K\}$, there are functions $w_k$ and $P_k$ satisfying
    \begin{equation*}
        (-1)^j u(j,R+k-j) = w_k(j) + P_k(j), \quad j \in I_k
    \end{equation*}
    where $P_k$ is a polynomial with $\deg P_k\leq k-1$, and
    \begin{equation*}
        \max_{1\leq k\leq K} \|w_k\|_{I_k} \leq M_u(D_R) \exp\left( CK\log\left( \frac{eR}{K} \right) \right).
    \end{equation*}
\end{lemma}
\begin{proof}
    For each $k\geq -1$ and each $j \in \mathbb{Z}$, let us consider
    \begin{equation*}
        v_k(j) = (-1)^j u(j,R+k-j).
    \end{equation*}
    Notice that $(j,R-1-j) \in D_R$ for all $j \in \{0,\dots,R-1\}$, which implies $|v_{-1}(j)| \leq M_u(D_R)$ for each $j \in \{0,\dots,R-1\}$. Similarly, we have that $|v_0(j)| \leq M_u(D_R)$ for every $j \in \{0,\dots,R\}$. Moreover, observe that
    \begin{equation*}
        \begin{aligned}
            (-1)^j u(j+1,R+k-1-j) &= -v_{k}(j+1),\\
            (-1)^j u(j-1,R+k-1-j) &= -v_{k-2}(j-1),\\
            (-1)^j u(j,R+k-j) &= v_{k}(j),\\
            (-1)^j u(j,R+k-2-j) &= v_{k-2}(j).\\
        \end{aligned}
    \end{equation*}
    Applying the harmonicity of $u$ at the point $(j,R+k-1-j)$ and multiplying by $(-1)^j$, we deduce
    \begin{equation*}
        4v_{k-1}(j) = -v_{k}(j+1) + v_k(j) -v_{k-2}(j-1) + v_{k-2}(j).
    \end{equation*}
    That is,
    \begin{equation}\label{eq:Remainder1}
        v_{k}(j+1) - v_k(j) = -4v_{k-1}(j)+ v_{k-2}(j)-v_{k-2}(j-1).
    \end{equation}
    Let us define
    \begin{equation*}
        w_{-1}(j) = \begin{cases}
            v_{-1}(j), &\text{if }0\leq j\leq R-1,\\
            0, &\text{otherwise},
        \end{cases} \quad \text{and}\quad w_{0}(j) = \begin{cases}
            v_{0}(j), &\text{if }0\leq j\leq R,\\
            0, &\text{otherwise}.
        \end{cases}
    \end{equation*}
    For each $k\geq 1$, let $w_k(-1) = w_k(0) = 0$, and let
    \begin{equation}\label{eq:RecursionW}
        w_k(j) = -4 \sum_{\ell=0}^{j-1}w_{k-1}(\ell) + w_{k-2}(j-1), \quad j\geq 1.
    \end{equation}
    Since $w_m(-1) = 0$ for each $m\geq -1$, the above formula also holds for $j=0$ if we interpret an empty sum as $0$. As a result,
    \begin{equation}\label{eq:Remainder2}
        \begin{aligned}
            w_k(j+1) - w_k(j) &= -4 \left(\sum_{\ell=0}^{j} w_{k-1}(\ell) - \sum_{\ell=0}^{j-1}w_{k-1}(\ell) \right) + w_{k-2}(j) - w_{k-2}(j-1)\\ &= -4w_{k-1}(j) + w_{k-2}(j) - w_{k-2}(j-1)
        \end{aligned}
    \end{equation}
    for every $j\geq 0$. Let us consider $r_k = v_k - w_k$. From \eqref{eq:Remainder1} and \eqref{eq:Remainder2}, we get
    \begin{equation}\label{eq:BaseCase1}
        r_{k}(j+1) - r_k(j) = -4r_{k-1}(j) + r_{k-2}(j) - r_{k-2}(j-1).
    \end{equation}
    Because of the choice of $w_{-1}$ and $w_0$, we deduce that 
    \begin{equation}\label{eq:BaseCase2}
        r_{-1}(j) = 0 \quad\text{for all }j \in \{0,\dots,R-1\}, \quad \text{and }r_{0}(j) = 0 \quad\text{for all }j \in \{0,\dots,R\}.
    \end{equation}
    Now, by using induction on $k \in \{1,\dots,K\}$, we shall prove that $r_k$ agrees with a polynomial of degree at most $k-1$ on $I_k$.

    For the base step, observe that \eqref{eq:BaseCase1} and \eqref{eq:BaseCase2} yield that $r_1(j+1) - r_1(j) = 0$ for all $j \in \{1,\dots,R-1\}$. Hence, $r_1(1) = \dots = r_1(R)$. Letting $P_1(t) = r_1(1)$, we get that $r_1(j) = P_1(j)$ for all $j \in I_1$ and $\deg P_1\leq 0$. This establishes the case $k=1$. 

    For the induction step, let $k\in\{2,\dots,K\}$, and let us assume that the desired polynomial representation holds for every $1\leq k'< k$. If $j\in \{k+1,\dots,R-2\}$ and $k\geq 3$, then $j\in I_k$ and $j,j-1 \in I_{k-1}$. By the induction hypothesis and \eqref{eq:BaseCase1}, we deduce
    \begin{equation*}
        r_k(j+1) - r_k(j) = -4P_{k-1}(j) + P_{k-2}(j) - P_{k-2}(j-1).
    \end{equation*}
    If $k=2$, the above equation also holds for every $j \in \{k+1,\dots,R-2\}$ by letting $P_0 \equiv 0$. That is,
    \begin{equation}\label{eq:IterationInduction}
        r_k(j+1) - r_k(j) = Q_k(j)
    \end{equation}
    for all $j \in \{k+1,\dots,R-2\}$, where $Q_k$ is the polynomial given by
    \begin{equation*}
        Q_k(t) = -4P_{k-1}(t) + P_{k-2}(t) - P_{k-2}(t-1).
    \end{equation*}
    Observe that the induction hypothesis implies that $\deg Q_k \leq k-2$. Then, there exists a polynomial $S_k$ with $\deg S_k \leq k-1$ such that $S_{k}(t+1) - S_k(t) = Q_k(t)$. Next, let $P_k(t) = S_k(t) + r_{k}(k+1) - S_k(k+1)$. This implies that $P_k(k+1) = r_k(k+1)$ and $\deg P_k\leq k-1$. By iterating \eqref{eq:IterationInduction}, we deduce
    \begin{equation*}
        \begin{aligned}
            r_k(j) &= r_k(k+1) + \sum_{\ell=k+1}^{j-1}Q_k(\ell)\\ &= P_k(k+1) + \sum_{\ell=k+1}^{j-1}(P_k(\ell + 1) - P_k(\ell))\\ &= P_k(j)
        \end{aligned}
    \end{equation*}
    for all $j \in I_k$ with $1\le k\le K$. This establishes the inductive step.

    Finally, we shall prove the desired bound for $w_k$. Let $a = 1 + K/R$ and $b = 8R/K$. By using induction on $k\geq -1$, we shall show that
    \begin{equation}\label{eq:InductionHypBoundW}
        |w_k(j)| \leq M_u(D_R) a^j b^{k+1}
    \end{equation}
    for all $j\geq 0$. Since $a>1$, we have
    \begin{equation*}
        |w_{-1}(j)| \leq  M_u(D_R) a^j
    \end{equation*}
    for all $j\geq 0$. Since $b \geq 16$, then
    \begin{equation*}
        |w_0(j)| \leq M_u(D_R)a^jb
    \end{equation*}
    for all $j\geq 0$. Now, let $k\geq 1$, and suppose that \eqref{eq:InductionHypBoundW} holds for each $-1\leq k' < k$. If $j=0$, the result follows immediately from $w_k(0)=0$. If $j\geq 1$, then \eqref{eq:RecursionW} and the induction hypothesis yield 
    \begin{equation*}
        \begin{aligned}
            |w_k(j)| &\leq 4M_u(D_R)b^{k}\sum_{\ell=0}^{j-1}a^{\ell} + M_u(D_R)a^{j-1}b^{k-1}\\ &\leq \frac{4M_u(D_R)a^jb^k}{a-1} + M_u(D_R)a^{j-1}b^{k-1}\\ &= M_u(D_R)a^{j}b^{k+1}\left( \frac{4}{b(a-1)} + \frac{1}{ab^2} \right)\\ &< M_u(D_R)a^jb^{k+1}.
        \end{aligned}
    \end{equation*}
    This establishes the inductive step.

    Let $k \in \{1,\dots,K\}$, and let $j \in I_k$. Notice that $|w_k(j)| \leq M_u(D_R) a^R b^{K+1}$, and
    \begin{equation*}
    \begin{aligned}
        \log(a^Rb^{K+1}) &= R \log\left( 1 + \frac{K}{R} \right) + (K+1)\log\left( \frac{8R}{K} \right)\\ &\leq K + 2K\log\left( \frac{8R}{K} \right)\\ &= 2K \log\left( \frac{eR}{K} \right) + (2\log 8-1)K\\ &\leq (2\log 8 +1)K\log\left( \frac{eR}{K} \right).
    \end{aligned}
    \end{equation*}
    Thus,
    \begin{equation*}
        |w_k(j)| \leq M_u(D_R) \exp\left( (2\log 8 +1) K \log\left( \frac{eR}{K} \right) \right).
    \end{equation*}
    This finishes the proof.
\end{proof}

We now focus our attention on estimating the polynomials $P_k$. Taking $k+1$ in \eqref{eq:BaseCase1} and using the polynomial representations on their respective intervals gives
\[
4P_k(j)
=
P_{k-1}(j)-P_{k-1}(j-1)
+
P_{k+1}(j)-P_{k+1}(j+1)
\]
for \(j\in\{k+2,\ldots,R-2\}\), with \(P_0=0\). Thus these polynomials satisfy the recurrence required in the next lemma.
\begin{lemma}\label{lemma:Step4}
    Let $K\geq 2$ and $R\geq 2K+1$ be integers. Let $P_0,\dots,P_K$ be polynomials with $P_0 = 0$ and $\deg P_k \leq k-1$ for each $k \in \{1,\dots,K\}$. Assume that, for each $k \in \{1,\dots,K-1\}$, the following identity holds
    \begin{equation}\label{eq:HypLemma4.4}
        4P_k(j) = P_{k-1}(j) - P_{k-1}(j-1) + P_{k+1}(j) - P_{k+1}(j+1)
    \end{equation}
    for every $j \in \{k+2,\dots,R-2\}$. Then,
    \begin{equation*}
        \max_{1\leq k < K} \sup_{t \in [0,R]} |P_k(t)| \leq \frac{K-1}{2} \sup_{t\in [-K,R+K]}|P_K(t)|.
    \end{equation*}
\end{lemma}

\begin{proof}
    For a polynomial $Q$, define
    \begin{equation*}
        D_+Q(t) = Q(t+1) - Q(t) \quad \text{and} \quad D_-Q(t) = Q(t) - Q(t-1).
    \end{equation*}
    We shall prove that \eqref{eq:HypLemma4.4} holds for every $t \in \mathbb{R}$. For each $k \in \{1,\dots,K-1\}$, let us consider
    \begin{equation*}
        A_k(t) = 4P_k(t) - D_-P_{k-1}(t) + D_+P_{k+1}(t).
    \end{equation*}
    Since $\deg P_{k+1}\leq k$, it follows that $\deg D_+P_{k+1} \leq k-1$. Therefore, $\deg A_k \leq k-1$. By \eqref{eq:HypLemma4.4}, we have that $A_k$ vanishes at the $R-k-3$ distinct points in the set $\{k+2,\dots,R-2\}$. Moreover, notice that the inequalities $R\geq 2K+1$ and $k\leq K-1$ imply that $R-k-3\geq k$. Given that $\deg A_k \leq k-1$ and $A_k$ has at least $k$ different zeros, we deduce that $A_k\equiv 0$. Hence,
    \begin{equation}\label{eq:RecurrencePolynomials}
        4P_k(t) = D_-P_{k-1}(t) - D_+P_{k+1}(t)
    \end{equation}
    for all $k \in \{1,\dots,K-1\}$.

    Now, let us consider the vector 
    \begin{equation*}
        \mathbf{P} = (P_1,\dots,P_{K-1}).
    \end{equation*}
    In addition, for a vector of polynomials $\mathbf{Z} = (Z_1,\dots,Z_{K-1})$, let us consider the transformation
    \begin{equation*}
        (T\mathbf{Z})_k = \frac{1}{4}(D_-Z_{k-1} - D_{+}Z_{k+1}), \quad k\in \{1,\dots,K-1\},
    \end{equation*}
    where $Z_0\equiv Z_K \equiv 0$. Let $\mathbf{F} = (0,\dots,0,- \frac{1}{4}D_{+}P_K)$. By \eqref{eq:RecurrencePolynomials}, we get that
    \begin{equation}\label{eq:RepPwithTP}
        \mathbf{P} = T\mathbf{P} + \mathbf{F}.
    \end{equation}
    Notice that $T$ lowers the maximum degree among the components by at least one. Since $\deg P_k\leq K-2$ for all $k \in \{1,\dots,K-1\}$, it follows that $T^{K-1}\mathbf{P} = 0$. By iterating \eqref{eq:RepPwithTP}, we deduce that
    \begin{equation}\label{eq:RepPwithTmF}
        \mathbf{P} = \sum_{m=0}^{K-2}T^m\mathbf{F}.
    \end{equation}
    Now, we aim to estimate each term in the above sum. For each $s \in \{0,\dots,K\}$, let $J_s = [-s,R+s]$, and define 
    \begin{equation*}
        \|\mathbf{Z}\|_s = \max_{1\leq k < K} \sup_{t \in J_s}|Z_k(t)|.
    \end{equation*}
    If $s\geq 1$ and $t \in J_{s-1}$, then $t-1,t,t+1 \in J_s$. Hence,
    \begin{equation*}
        \sup_{t \in J_{s-1}}|(T\mathbf{Z})_k(t)| \leq \frac{1}{2}\left( \sup_{t\in J_s}|Z_{k-1}(t)| + \sup_{t\in J_s}|Z_{k+1}(t)| \right) \leq \|\mathbf{Z}\|_s
    \end{equation*}
    for each $k \in \{1,\dots,K-1\}$. Therefore,
    \begin{equation*}
        \|T\mathbf{Z}\|_{s-1} \leq \|\mathbf{Z}\|_s.
    \end{equation*}
    If $t \in J_{K-1}$, then $t,t+1 \in J_K$. As a result,
    \begin{equation*}
        \|\mathbf{F}\|_{K-1} = \sup_{t \in J_{K-1}}|F_{K-1}(t)| \leq \frac{1}{2} \sup_{t\in J_K}|P_K(t)|.
    \end{equation*}
    Combining the last two inequalities, we get
    \begin{equation*}
        \|T^m\mathbf{F}\|_{K-1-m} \leq \|\mathbf{F}\|_{K-1} \leq \frac{1}{2} \sup_{t\in J_K}|P_K(t)|,
    \end{equation*}
    for each $m \in \{0,\dots,K-2\}$. Since the interval $J_{K-1-m}$ contains $J_0 = [0,R]$ for every $m \in \{0,\dots,K-2\}$, an application of the triangle inequality in \eqref{eq:RepPwithTmF} yields
    \begin{equation*}
        \|\mathbf{P}\|_0 \leq \sum_{m=0}^{K-2} \|T^m\mathbf{F}\|_0 \leq \sum_{m=0}^{K-2}\|T^m\mathbf{F}\|_{K-1-m} \leq \frac{K-1}{2}\sup_{t\in J_K}|P_K(t)|. 
    \end{equation*}
    This finishes the proof.
\end{proof}

Therefore, it remains to estimate $P_K$. For this step, the non-negativity assumption on $u$ will play a fundamental role. Indeed, we show that if $u$ takes negative values at a small number of points on the line $(j,R+K-j)$ with $j \in I_K$, then we can deduce the desired bound for $P_K$.

\begin{lemma}\label{lemma:Step3}
    Let $R\geq 4$ and $1\leq K\leq R/2-1$ be integers, and let $P_K$ be the polynomial from Lemma \ref{Lemma:4.2}. Assume that
    \begin{equation}\label{eq:HypLemma4.3}
        t_K := |\{j\in I_K : u(j,R+K-j) < 0\}| \leq \frac{R-2K-2}{2}.
    \end{equation}
    Then,
    \begin{equation*}
        \sup_{t \in [-K,R+K]}|P_K(t)| \leq M_u(D_R) \exp\left( C K  \log\left( \frac{eR}{K} \right)\right),
    \end{equation*}
    where $C>0$ is a constant independent of $R$ and $K$.
\end{lemma}

\begin{proof}
    Notice that $L := |I_K| = R-K-1$. Let 
    \begin{equation*}
        I_K \setminus \{j\in I_K : u(j,R+K-j) < 0\} = \{a_1,\dots,a_s\},
    \end{equation*}
    where $s = L - t_K$ and $a_1 < \dots<a_s$. As a result, $u(a_i,R+K-a_i) \geq 0$ for each $i \in \{1,\dots,s\}$. Observe that if $a_i$ and $a_{i+1}$ have the same parity, there is at least one integer between them. Hence, there are at most $t_K$ pairs $(a_i,a_{i+1})$ with the same parity. Moreover, notice that the number of maximal runs of consecutive entries of the ordered sequence \(\{a_1,\ldots,a_s\}\) having the same parity is equal to $s$ minus the number of consecutive pairs $(a_i,a_{i+1})$ with the same parity, which is at least $s - t_K$. By \eqref{eq:HypLemma4.3}, we deduce
    \begin{equation*}
        s-t_K = L-2t_K \geq R-K-1 - (R-2K-2) = K+1.
    \end{equation*}
    Let us choose one integer from each of the first $K+1$ maximal blocks of consecutive entries of the ordered sequence $a_1,\ldots,a_s$ having the same parity, say $j_0 < \dots < j_K$. Thus,
    \begin{equation*}
        u(j_i,R+K-j_i)\geq 0 \quad\text{and}\quad j_{i} \equiv j_0+i \mod 2
    \end{equation*}
    for all $i \in \{0,\dots,K\}$. For each $i \in \{0,\dots,K\}$, let us define $q_i = (-1)^{j_i}P_K(j_i)$. By Lemma \ref{Lemma:4.2}, we get
    \begin{equation*}
        q_i = u(j_i,R+K-j_i) - (-1)^{j_i} w_K(j_i) \geq -M_u(D_R) \exp\left( CK \log\left( \frac{eR}{K} \right) \right)
    \end{equation*}
    for each $i \in \{0,\dots,K\}$. Moreover, Lagrange interpolation yields
    \begin{equation*}
        P_K(t) = \sum_{i=0}^{K} P_K(j_i) \prod_{\substack{\ell=0\\ \ell\neq i}}^{K} \frac{t-j_{\ell}}{j_i-j_{\ell}},
    \end{equation*}
    since the polynomial on the right-hand side has degree at most $K$ and agrees with $P_K$ at $K+1$ different points. Since $\deg P_K\leq K-1$, it follows that the coefficient corresponding to $t^K$ in the polynomial on the right-hand side is zero, i.e.,
    \begin{equation}\label{eq:SumNodes0}
        \sum_{i=0}^{K} P_K(j_i) \prod_{\substack{\ell=0\\ \ell\neq i}}^{K} \frac{1}{j_i-j_{\ell}} = 0.
    \end{equation}
    For each $i \in \{0,\dots,K\}$, let us consider
    \begin{equation*}
        \omega_i =  \prod_{\substack{\ell=0\\ \ell\neq i}}^{K} \frac{1}{|j_i-j_{\ell}|}.
    \end{equation*}
    Notice that there are exactly $K-i$ negative numbers in the product $ \prod_{\substack{\ell=0\\ \ell\neq i}}^{K} (j_i-j_{\ell})$, which implies
    \begin{equation*}
        P_K(j_i) \prod_{\substack{\ell=0\\ \ell\neq i}}^{K} \frac{1}{j_i-j_{\ell}} = (-1)^{K-i}P_K(j_i) \omega_i = (-1)^{K-i+j_i}\omega_iq_i = (-1)^{K+j_0}\omega_i q_i,
    \end{equation*}
    where in the last equation we used the fact that $j_i-i$ and $j_0$ have the same parity. Applying \eqref{eq:SumNodes0}, we get
    \begin{equation*}
        \sum_{i=0}^K \omega_i q_i = 0.
    \end{equation*}
    Let $i \in \{0,\dots,K\}$ be fixed. Then,
    \begin{equation*}
        \omega_i q_i = - \sum_{\substack{\ell=0\\ \ell\neq i}}^{K}\omega_{\ell}q_{\ell} \leq M_u(D_R) \exp\left( CK \log\left( \frac{eR}{K} \right) \right) \sum_{\substack{\ell=0\\ \ell\neq i}}^{K}\omega_{\ell}.
    \end{equation*}
    As a result,
    \begin{equation}\label{eq:BoundPk1}
        |P_K(j_i)| = |q_i| \leq M_u(D_R) \exp\left( CK \log\left( \frac{eR}{K} \right) \right) \frac{\sum_{\ell=0}^{K}\omega_{\ell}}{\omega_i}.
    \end{equation}
    Since $|I_K| \leq R$ and all the $j_\ell$'s belong to $I_K$, it follows that 
    \begin{equation}\label{eq:BoundWeight1}
        \frac{1}{\omega_i} = \prod_{\substack{\ell=0\\ \ell\neq i}}^{K} |j_i-j_{\ell}| \leq R^K.
    \end{equation}
    Additionally, since the integers $j_\ell$ are chosen in increasing order, we have that $j_i-j_{\ell} \geq i-\ell$ for all $\ell < i$, and $j_{\ell} - j_i \geq \ell-i$ for all $\ell>i$. Therefore,
    \begin{equation}\label{eq:Spacing}
        \prod_{\substack{\ell=0\\ \ell\neq i}}^{K} |j_i-j_{\ell}| \geq i! (K-i)!.
    \end{equation}
    Thus,
    \begin{equation}\label{eq:BoundWeight2}
        \sum_{i=0}^K \omega_i \leq \sum_{i=0}^K \frac{1}{i!(K-i)!} = \frac{1}{K!}\sum_{i=0}^K \binom{K}{i} = \frac{2^K}{K!}.
    \end{equation}
    Plugging \eqref{eq:BoundWeight1} and \eqref{eq:BoundWeight2} into \eqref{eq:BoundPk1}, we deduce
    \begin{equation}\label{eq:BoundPkNodes}
        |P_K(j_i)| \leq M_u(D_R) \exp\left( CK \log\left( \frac{eR}{K} \right) \right) \frac{(2R)^K}{K!}.
    \end{equation}
    Now, we shall extend this bound beyond the nodes $j_i$. Since $\deg P_K\leq K-1$, an application of Lagrange interpolation yields
    \begin{equation}\label{eq:LagrangeRepPk}
        P_K(t) = \sum_{i=0}^{K-1} P_K(j_i) \prod_{\substack{\ell=0\\ \ell\neq i}}^{K-1} \frac{t-j_{\ell}}{j_i-j_{\ell}}.
    \end{equation}
    If $t \in [-K,R+K]$, the fact that $j_{\ell} \in [0,R]$ implies that $|t-j_\ell| \leq R+K \leq 2R$. Arguing as in \eqref{eq:Spacing}, we can also deduce that
    \begin{equation*}
        \prod_{\substack{\ell=0\\ \ell\neq i}}^{K-1} |j_i-j_{\ell}| \geq i! (K-1-i)!.
    \end{equation*}
    Applying these bounds and \eqref{eq:BoundPkNodes} in \eqref{eq:LagrangeRepPk}, we have
    \begin{equation}\label{eq:FinalPk1}
        \begin{aligned}
            |P_K(t)| &\leq M_u(D_R) \exp\left( CK \log\left( \frac{eR}{K} \right) \right) \frac{(2R)^K}{K!} (2R)^{K-1}\sum_{i=0}^{K-1} \frac{1}{i! (K-1-i)!}\\ &=M_u(D_R) \exp\left( CK \log\left( \frac{eR}{K} \right) \right) \frac{(2R)^K}{K!} (2R)^{K-1} \frac{2^{K-1}}{(K-1)!}\\ &= M_u(D_R) \exp\left( CK \log\left( \frac{eR}{K} \right) \right) \frac{(2R)^K(4R)^{K-1}}{K!(K-1)!}.
        \end{aligned}
    \end{equation}
    Finally, notice that the inequality $K! \geq (K/e)^K$ implies
    \begin{equation}\label{eq:FinalPk2}
        \begin{aligned}
            \frac{(2R)^K(4R)^{K-1}}{K!(K-1)!} &= \frac{K}{4R}\frac{(8R^2)^K}{(K!)^2}\\ 
            &\leq \left( \frac{8e^2R^2}{K^2} \right)^K\\ &= \exp\left( K\left( \log 8 + 2 \log\left( \frac{eR}{K} \right) \right) \right)\\ &\leq \exp\left( 5K \log\left( \frac{eR}{K} \right)\right).
        \end{aligned}
    \end{equation}
    Combining \eqref{eq:FinalPk1} and \eqref{eq:FinalPk2}, we obtain the desired result.
\end{proof}

It remains to choose $K$ so that the estimate for $P_K$ applies in all four quadrants. The next lemma uses the upper density of $\{u<0\}$ to show that, for every sufficiently large $R$, there exists $K\in[\gamma R,(1/2-\gamma)R]$ satisfying \eqref{eq:HypLemma4.3} on the four outer segments. This gives the required bound on the outer polynomial in each quadrant. The polynomial recurrence then controls the intermediate diagonals and implies Proposition \ref{prop:large-range}.

\begin{lemma} \label{lemma:Step1}
    For each $\varepsilon = (\varepsilon_1 , \varepsilon_2) \in \{-1,1\}^2$, define
    \begin{equation}\label{def:SetSrk}
        S_{R,K}^{\varepsilon} = \{(\varepsilon_1 j , \varepsilon_2(R+K-j)) : j \in I_K\}.
    \end{equation}
    Then, there exists $R_0 = R_0(u) \in \mathbb{Z}^+$ such that for every integer $R\geq R_0$ there exists $K \in \mathbb{Z}^+$ with
    \begin{equation*}
        \gamma R \leq K \leq \left( \frac{1}{2} - \gamma \right)R
    \end{equation*}
    satisfying
    \begin{equation*}
        |S_{R,K}^{\varepsilon} \cap \{u<0\}| \leq \frac{R-2K-2}{2}
    \end{equation*}
    for every $\varepsilon \in \{-1,1\}^2$.
\end{lemma}

\begin{proof}
    Let $\widetilde{R}_0\geq 8$ be a sufficiently large integer such that 
    \begin{equation*}
        \left[ \gamma \widetilde{R}_0 , \left( \frac{1}{2} - \gamma \right)\widetilde{R}_0 \right] \cap \mathbb{Z} \neq \emptyset.
    \end{equation*}
    Let $\mathcal{K}$ be the set of integers $R\geq \widetilde{R}_0$ such that for each
    \begin{equation*}
        K \in \left[ \gamma R , \left( \frac{1}{2} - \gamma \right)R \right] \cap \mathbb{Z}
    \end{equation*}
    there exists $\varepsilon(K) \in \{-1,1\}^2$ satisfying
    \begin{equation*}
        |S_{R,K}^{\varepsilon (K)} \cap \{u<0\}| > \frac{R-2K-2}{2} = \frac{R}{2} - K - 1.
    \end{equation*}
    For the sake of contradiction, let us assume that $\mathcal{K}$ is infinite. Let $R \in \mathcal{K}$ be fixed, and let 
    \begin{equation*}
        a_R = \lceil \gamma R \rceil, \quad \text{and} \quad b_R =  \left\lfloor\left( \frac{1}{2} - \gamma \right)R \right\rfloor.
    \end{equation*}
    Since the segments $S_{R,K}^{\varepsilon (K)} \cap \{u<0\}$ are contained in $Q_R$ and are pairwise disjoint (see Figure \ref{fig1}), we get
    \begin{equation}\label{eq:LowerBoundDeltaU}
        |Q_R\cap \{u<0\}| \geq \sum_{K=a_R}^{b_R}|S_{R,K}^{\varepsilon (K)} \cap \{u<0\}| > \sum_{K=a_R}^{b_R}\left( \frac{R}{2} - K - 1 \right).
    \end{equation}
    Let $m_R = |[a_R,b_R] \cap \mathbb{Z}| = b_R - a_R + 1$. Notice that 
    \begin{equation*}
        a_R + b_R \leq \frac{R}{2} + 1 \quad \text{and} \quad m_R \geq \left( \frac{1}{2} - 2\gamma \right) R - 1.
    \end{equation*} 
    Thus,
    \begin{equation}\label{eq:LowerBoundDeltaU2}
    \begin{aligned}
        \sum_{K=a_R}^{b_R}\left( \frac{R}{2} - K - 1 \right) &= m_R \left( \frac{R}{2} - 1 - \frac{a_R + b_R}{2}\right)\\ &\geq \left[\left( \frac{1}{2} - 2\gamma \right) R - 1 \right] \left( \frac{R}{4} - \frac{3}{2} \right)\\ &= \left( \frac{1}{8} - \frac{\gamma}{2} \right) R^2 - \left( 1-3\gamma \right)R + \frac{3}{2}\\ &\geq \left( \frac{1}{8} - \frac{\gamma}{2} \right) R^2 - R.
    \end{aligned}
    \end{equation}
    From \eqref{eq:LowerBoundDeltaU} and \eqref{eq:LowerBoundDeltaU2}, we deduce
    \begin{equation*}
        \delta_R^-(u) = \frac{|Q_R \cap \{u<0\}|}{(2R+1)^2} \geq \frac{\left( \frac{1}{8} - \frac{\gamma}{2} \right) R^2 - R}{(2R+1)^2}
    \end{equation*}
    for all $R \in \mathcal{K}$. Since $\mathcal{K}$ is infinite, there exists an increasing sequence $\{R_k\}_{k\geq 1} \subseteq \mathcal{K}$ such that $R_k \to \infty$ as $k \to \infty$. It follows that
    \begin{equation*}
        \delta = \overline{\delta}^-(u) \geq \liminf_{k \to \infty}\delta_{R_k}^-(u) \geq \frac{1}{32} - \frac{\gamma}{8} = \frac{1}{64} + \frac{\delta}{2} > \delta,
    \end{equation*}
    which is a contradiction. We finish the proof by letting $R_0 = \max(\mathcal{K}\cup \{\widetilde{R}_0\}) + 1$.
\end{proof}

We can now deduce Proposition \ref{prop:large-range} as a simple consequence of the previous lemmas.

\begin{proof}[Proof of Proposition \ref{prop:large-range}]
    Let $u : \mathbb{Z}^2 \to \mathbb{R}$ be a harmonic function with $\delta = \overline{\delta}^-(u) < 1/32$. Thus, Lemma \ref{lemma:Step1} yields the existence of $R_0 = R_0(u) \in \mathbb{Z}^+$ such that for every integer $R\geq R_0$ there exists $K \in \mathbb{Z}^+$ with
    \begin{equation*}
        \gamma R \leq K \leq \left( \frac{1}{2} - \gamma \right)R
    \end{equation*}
    satisfying
    \begin{equation*}
        |S_{R,K}^{\varepsilon} \cap \{u<0\}| \leq \frac{R-2K-2}{2}
    \end{equation*}
    for every $\varepsilon \in \{-1,1\}^2$ (where $S_{R,K}^{\varepsilon}$ is defined as in \eqref{def:SetSrk}). Let $Q_I = \{(x,y) \in \mathbb{Z}^2 : x,y\geq 0\}$. Without loss of generality, we can assume that $R_0 \geq 2/\gamma$. As a result, if $R\geq R_0$ and $\gamma R\leq K\leq (1/2-\gamma)R$, then $2\leq K\leq R/2-2$. Using the decomposition from Lemma \ref{Lemma:4.2} together with the bounds from Lemmas \ref{lemma:Step3} and \ref{lemma:Step4}, we deduce
    \begin{equation*}
        \|u\|_{D_{R+K} \cap Q_{R-1} \cap Q_{I}} \leq M_u(D_R) \exp\left( CK\log\left( \frac{eR}{K} \right) \right)
    \end{equation*}
    for some constant $C$ independent of $u$ and $\delta$. By composing $u$ with rotations of $\pi/2$, $\pi$, and $3\pi/2$ (as in \eqref{eq:RotationsU}), we can combine Lemmas \ref{Lemma:4.2}, \ref{lemma:Step3} and \ref{lemma:Step4} as before to establish analogous bounds for $u$ in the other three quadrants of the plane. Hence,
    \begin{equation*}
        \|u\|_{D_{R+K} \cap Q_{R-1}} \leq M_u(D_R) \exp\left( CK\log\left( \frac{eR}{K} \right) \right).
    \end{equation*}
    Notice that this constant $C$ is not necessarily the same as before, but it is still independent of $u$ and $\delta$. Therefore, Lemma \ref{lemma:Step5} implies that
    \begin{equation}\label{eq:AbsorbPowerIntoExp}
    \begin{aligned}
        \|u\|_{D_{R+K}} &\leq 7^{K+1} M_u(D_R) \exp\left( CK\log\left( \frac{eR}{K} \right) \right)\\ &\leq M_u(D_R)\exp\left( (C + 2\log(7)) K \log\left( \frac{eR}{K} \right)\right).
    \end{aligned}
    \end{equation}
    Since $M_u(D_R)\geq1$ and the exponential factor is at least $1$, the preceding estimate also holds for $\displaystyle M_u(D_{R+K})= \max\left\{1,\|u\|_{D_{R+K}}\right\}.$  Applying Lemma \ref{lemma:InitialReduction} completes the proof.
\end{proof}
\subsection{Proof of Theorem \ref{thm:main}}
For each $z\in\Z^2$, let us consider the translated diamonds and squares
\begin{equation*}
    D_r(z)=z+D_r \quad \text{and}\quad  Q_r(z)=z+Q_r.
\end{equation*}
Moreover, let $\mathcal N=\{x \in \mathbb{Z}^2:u(x)<0\}$, $\delta_r=\delta_r^-(u)$ and $M_u(D_r(z))=\max\{1,\|u\|_{D_r(z)}\}$ for every $r\in \mathbb{N}$. The next lemma shows that if $u$ takes negative values in a small portion of a diagonal in a translated diamond, then we can estimate the maximum of $|u|$ in that diamond by an exponential factor of its maximum in a smaller diamond, which is analogous to \eqref{eq:AbsorbPowerIntoExp}. For any $2\leq K\leq R/2-2$ and $\varepsilon = (\varepsilon_1 , \varepsilon_2) \in \{-1,1\}^2$, we define translated diagonal
\begin{equation*}
    S_{R,K}^{\varepsilon}(z)= z + S_{R,K}^{\varepsilon},
\end{equation*}
where $S_{R,K}^{\varepsilon}$ is the diagonal defined in \eqref{def:SetSrk}.

\begin{lemma}\label{lem:local-diamonds} Let $z\in\mathbb Z^2$, and let $R,K$ be integers with $2\le K\le R/2-2$. Suppose that 
\begin{equation}
 |S_{R,K}^{\sigma}(z)\cap\mathcal N|\le R/2-K-1 
 \label{sign}
\end{equation}
for every $\sigma\in\{-1,1\}^2$. Then,
\begin{equation}
 M_u(D_{R+K}(z))\leq M_u(D_R(z))\exp\left(CK\log\frac{eR}{K}\right). \label{eq:local-diamonds}
\end{equation}
\end{lemma}
\begin{proof}
Let $\sigma\in\{-1,1\}^2$ be fixed, and let us consider $v_\sigma(a,b)=u(z+(\sigma_1a,\sigma_2b))$. Notice that $v_{\sigma}$ is also harmonic and $M_{v_{\sigma}}(D_R) = M_u(D_R(z))$ for every $R$. The desired result can be deduced by combining Lemmas \ref{Lemma:4.2}, \ref{lemma:Step3}, \ref{lemma:Step4}, and \ref{lemma:Step5} as in the proof of Proposition \ref{prop:large-range}.
\end{proof}

Next, we use this exponential growth estimate on translated concentric diamonds that lie inside squares with small negative density. More precisely, fixing $R_*=64$ and choosing $\eta>0$ with $9\eta<1/64$ and $\eta(2R_*+3)^2<1$, we define the following collection of centers
\begin{equation}\label{def:centers_small_neg}
  G:=\{z \in \mathbb{Z}^2:|\mathcal N\cap Q_S(z)|\leq\eta(2S+1)^2\text{ for every }S\geq 0\}.   
\end{equation}
Arguing as in the proof of Lemma~\ref{lemma:Step1}, we may prove that for every $z\in G$ and every $R\geq R_*$ there is an integer $K\in[R/8,R/4]$ satisfying the sign condition \eqref{sign} in all four quadrants. 

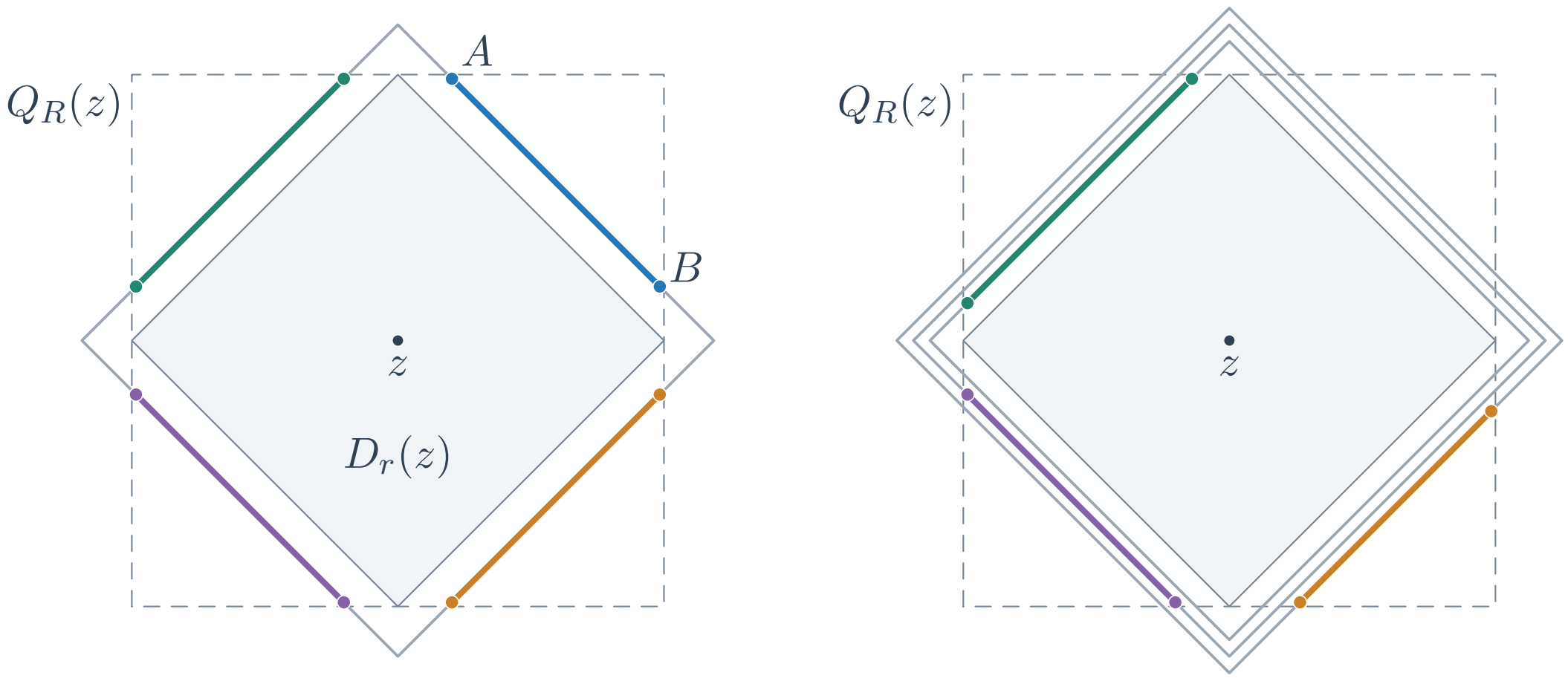
\begin{figure}[htbp]
\centering
\begin{tikzpicture}[
  x=.037cm,y=.037cm,font=\small,text=LabelColor,
  line cap=round,line join=round,
  square/.style={draw=Boundary,dash pattern=on 4pt off 3pt},
  outer/.style={draw=OuterBoundary,line width=.7pt}]
\def\diamond#1{(#1,0)--(0,#1)--(-#1,0)--(0,-#1)--cycle}
\def\basepanel{%
  \filldraw[fill=InnerFill,draw=Boundary] \diamond{64};
  \draw[square] (-64,-64) rectangle (64,64);
  \fill[LabelColor] (0,0) circle[radius=1.3pt];
  \node[below] at (0,0) {$z$};}
\def\segment#1#2#3#4{%
  \begin{scope}[xscale=#2,yscale=#3]
    \draw[#4,line width=1.5pt] ({#1+1},63)--(63,{#1+1});
    \foreach \j in {#1+1,63}{
      \filldraw[fill=#4,draw=white,line width=.3pt]
        ({\j},{64+#1-(\j)}) circle[radius=1.7pt];}
  \end{scope}}

\begin{scope}
  \basepanel
  \draw[outer] \diamond{76};
  \segment{12}{1}{1}{NEblue}
  \segment{12}{-1}{1}{NWgreen}
  \segment{12}{-1}{-1}{SWpurple}
  \segment{12}{1}{-1}{SEorange}
  \node[above right,inner sep=1pt] at (14,65) {$A$};
  \node[above right,inner sep=1pt] at (64,13) {$B$};
  \node at (0,-28) {$D_r(z)$};
  \node[anchor=north east,inner sep=2pt] at (-64,64) {$Q_R(z)$};
  \node at (0,-91) {(a) One candidate};
\end{scope}
\begin{scope}[shift={(200,0)}]
  \basepanel
  \foreach \t in {72,76,80}{\draw[outer] \diamond{\t};}
  \segment{8}{-1}{1}{NWgreen}
  \segment{12}{-1}{-1}{SWpurple}
  \segment{16}{1}{-1}{SEorange}
  \node at (0,-91) {(b) Distinct candidates};
    \node[anchor=north east,inner sep=2pt] at (-64,64) {$Q_R(z)$};
\end{scope}
\end{tikzpicture}
\caption{Selection of a translated diamond, illustrated with $R=64$. In (a), the four colored segments are tested for the same $K=12$. The endpoints of the blue edge are $A=z+(K+1,R-1)$ and $B=z+(R-1,K+1)$. In (b), suppose that the three segments fail the sign condition for $K=8,12,16$, respectively. They lie in $Q_R(z)$ at distinct $\ell^1$-distances from $z$, and thus their negative lattice points can be counted separately as in \eqref{eq:LowerBoundDeltaU}-\eqref{eq:LowerBoundDeltaU2}.}
\label{fig:translated-diamond-selection}
\end{figure}

The next lemma gives us a precise exponential growth estimate for $|u|$ around the chosen centers.

\begin{lemma}\label{lem:good-center-growth}
Let $g(x)=\log(1+|u(x)|)$. There is an absolute constant $A$ such
that, for every $z\in G$ and every $w\in\mathbb Z^2$,
\begin{equation}\label{cone}
 g(w)\leq g(z)+A|w-z|_1.
\end{equation}
\end{lemma}

\begin{proof}
Let $z\in G$ and $R_0=R_*$. By the preceding counting argument, at each radius $R_j$ we may choose an integer $K_j\in[R_j/8,R_j/4]$ satisfying \eqref{sign} in all four quadrants. Applying Lemma \ref{lem:local-diamonds} and the iteration from the proof of Lemma~\ref{lemma:InitialReduction} to the translated harmonic function $v_\sigma$ and the sequence $R_{j+1}=R_j+K_j$, we obtain
\begin{equation}\label{eq:diamond-growth-iterative}
 M_u(D_S(z))\leq M_u(D_{R_0}(z))e^{C_1S}
\end{equation}
for every $S \in \mathbb{Z}^+$. 

Now we shall connect the value $M_u(D_{R_0}(z))=M_u(D_{R_*}(z))$ from the initial diamond to the magnitude of $u$ at its center. Since $z\in G$, the definition of $G$ yields
\begin{equation*}
    |\mathcal N\cap Q_{R_0+1}(z)|
 \leq\eta(2R_0+3)^2<1.
\end{equation*}
Notice that left-hand side of the above equation is a non-negative integer, which means that it equals zero. Hence, $u(x,y)\geq 0$ for all $(x,y) \in Q_{R_0+1}(z)$. Let $w\in D_{R_0}(z)$, and
let us consider a shortest lattice path from $z$ to $w$: $z=x_0,x_1,\ldots,x_{|w-z|_1}=w$. Observe that each vertex of the path lies in $D_{R_0}(z)$, and all four of its neighbors lie in $Q_{R_0+1}(z)$. For each successive pair $x,x'$ on the path, the harmonicity and non-negativity of $u$ yield
\begin{equation*}
    u(x')\leq\sum_{|y-x|_1=1}u(y)=4u(x).
\end{equation*}
It follows that $0\leq u(w)\leq 4^{|w-z|_1}u(z)$. Since $u(z)\geq 0$, we have that $g(z)=\log(1+u(z))$. Hence,
\begin{equation*}
    1+|u(w)|\leq 1+4^{|w-z|_1}u(z)\leq 4^{|w-z|_1}(1+u(z)),
\end{equation*}
and thus $g(w)\leq g(z)+\log4|w-z|_1$. On the other hand, if $|w-z|_1 > R_{\ast}$, then \eqref{eq:diamond-growth-iterative} yields
\begin{equation*}
    1+|u(w)|\leq 2M_u(D_{|w-z|_1}(z)) \leq 2\cdot4^{R_*}(1+u(z))e^{C_1|w-z|_1}.
\end{equation*}
Taking logarithms and using that $|w-z|_1 > R_*\geq1$, we get
\begin{equation*}
    g(w)-g(z) \leq \log2+R_*\log4+C_1|w-z|_1 \leq (\log2+\log4+C_1)|w-z|_1 ,
\end{equation*}
which proves \eqref{cone} with absolute constant $A=C_1+\log 8$.
\end{proof}

Now that we have exponential growth control for the centers in $G$, we count the number of the remaining centers in $Q_R\setminus G$. For every sufficiently large $R$, we show that if $z\in Q_R\setminus G$, then any square for which \eqref{def:centers_small_neg} fails is contained in $Q_{2R}$.

Let $\delta_U=\eta/900$, and let us suppose that $\overline{\delta}^-(u)\leq\delta_U$. Thus, there exists
$R_u$ such that $\delta_T\leq 2\delta_U$ for every $T\geq R_u$.
Let $R\geq R_u$. If $z\in Q_R$ and $S\geq R$, then
$Q_S(z)\subseteq Q_{S+R}\subseteq Q_{2S}$. As a result,
\begin{equation}\label{weak1}
 \frac{|\mathcal N\cap Q_S(z)|}{|Q_S(z)|}
 \leq \delta_{2S}\frac{|Q_{2S}|}{|Q_S|}
 \leq 2\delta_U\frac{(4S+1)^2}{(2S+1)^2}
 <10\delta_U.    
\end{equation}
For each $z\in Q_R\setminus G$, we can choose a radius $S_z$ such that $|\mathcal N\cap Q_{S_z}(z)|>\eta|Q_{S_z}(z)|$. Notice that the previous estimate implies that $S_z<R$, and thus $Q_{S_z}(z) \subseteq Q_{2R}$. Since $\{Q_{S_z}(z)\}_{z \in Q_R\setminus G}$ is a finite cover of $Q_R\setminus G$, an application of Vitali's covering lemma yields the existence of a $w_1,\dots,w_J \in Q_R\setminus G$ such that the squares $\{Q_{S_{w_i}}(w_i)\}_{1\leq i\leq J}$ are pairwise disjoint and satisfy 
\begin{equation*}
    |Q_R\setminus G| \leq  9\sum_{i=1}^{J} |Q_{S_{w_i}}(w_i)|.
\end{equation*}
By the choice of the radii $S_{w_i}$ and the pairwise disjointness, we get
\begin{equation}\label{weak}
    |Q_R\setminus G| \leq \frac{9}{\eta}\sum_{i=1}^{J}|\mathcal N\cap Q_{S_{w_i}}(w_i)| \leq \frac{9}{\eta}\,|\mathcal N\cap Q_{2R}|.
\end{equation}

Dividing by $|Q_R|$ and using \eqref{weak1}, we obtain
\[
 \frac{|Q_R\setminus G|}{|Q_R|}
 \leq \frac{45}{\eta}\,\frac{|\mathcal N\cap Q_{2R}|}{|Q_{2R}|}
 < \frac{450\delta_U}{\eta}=\frac{1}{2}.
\]
Thus, $|G\cap Q_R|\geq |Q_R|/2$ for every $R\ge R_u$. In particular, $G\neq \emptyset$. This completes the count of exceptional centers.

The following result constructs a global \textit{Lipschitz majorant} that agrees with $g$ on the good centers and whose cubic energy is controlled by the number of negative sites.

\begin{lemma}\label{lem:lipschitz-energy}
Let $\delta\le\delta_U$ and $R_u$ be as in \eqref{weak1}, let $g$ be the function defined in Lemma \ref{lem:good-center-growth}. Define
\begin{equation*}
    F(x):=\inf_{z\in G}\bigl(g(z)+A|x-z|_1\bigr)
\end{equation*}
for each $x \in \mathbb{Z}^2$, where $A$ is the constant in \eqref{cone}. Then, $F$ is finite and $A$-Lipschitz on $\mathbb Z^2$, with $F\geq g$ everywhere and $F=g$ on $G$.

If $\mathcal N=\emptyset$, then $u$ and $F$ are constant.
If $\mathcal N\neq\emptyset$ and
$d(x)=\operatorname{dist}_\infty(x,\mathcal N)$, then for every pair of neighboring points $x,y\in G$,
\begin{equation}\label{dist}
 |F(x)-F(y)|\le \frac{C}{1+d(x)}.
\end{equation}
In either case, for every integer $R\ge\max\{1,R_u\}$,
\begin{equation}\label{energy}
 \sum_{x\in Q_{2R}}\sum_{i=1}^2 |F(x+e_i)-F(x)|^3
 \leq C\left(|\mathcal N\cap Q_{6R}|+\frac1R\right).
\end{equation}
\end{lemma}
\begin{proof}
We start verifying the claimed properties of $F$ on $\mathbb Z^2$. By the preceding covering estimate, we know that $G\ne\emptyset$. Let $z_0\in G$.
Taking the infimum in \eqref{cone} gives
\[
 g(x)\le F(x)\le g(z_0)+A|x-z_0|_1<\infty.
\]
Thus, $F$ is finite and $g\leq F$. If $x\in G$, letting $z=x$ in the definition of $F$ gives the reverse inequality $F(x)\le g(x)$, which means that $F$ agrees with $g$ on $G$. Moreover, let $x,y \in \mathbb{Z}^2$. Notice that for every $z\in G$,
\[
 F(x)\leq g(z)+A|x-z|_1\le g(z)+A|y-z|_1+A|x-y|_1,
\]
which implies that $F(x) \leq F(y) + A|x-y|_1$. Interchanging $x$ and $y$, we get $|F(x)-F(y)|\leq A|x-y|_1$.

Let us consider $h=1+u$. The relative interior difference estimate
\cite[Theorem~6.3.8(b), (6.19)]{LawlerLimicBook} says that if $h$ is nonnegative on a lattice ball of sufficiently large radius $n$ and on its lattice boundary, and if $h$ is harmonic inside the ball, then
\[
 |h(q)-h(p)|\le \frac{C}{n}h(p)
\]
for its center $p$ and every neighbor $q$ of $p$.
If $\mathcal N=\emptyset$, then $h\ge1$ on the whole $\mathbb{Z}^2$. Letting $n\to\infty$, we get that $h$ is constant. Therefore, $u$ and $F$ are constant. 

Suppose now that $\mathcal N\ne\emptyset$, and let $x,y\in G$ be neighbors. Since $F=g$ on $G$ and $u\ge0$ there, we have $F(x)=\log h(x)$ and $F(y)=\log h(y)$.  For $d(x)$ sufficiently large, let us consider $n=\lfloor d(x)/3\rfloor$. The ball of radius $n$ centered at $x$, together with its lattice boundary, lies in $Q_{n+1}(x)$. Since $n+1<d(x)$, it follows that this ball has empty intersection with $\mathcal{N}$, which implies that $h\geq 1$ throughout it. Applying the relative difference estimate at $x$ gives
\[
 |h(x)-h(y)|\leq\frac{C_0}{d(x)}h(x).
\] 
For $d(x)$ large enough such that $C_0/d(x)\le 1/2$, the relative difference above implies $h(y)\ge h(x)/2$. As a result,
\[
 |h(x)-h(y)|
 \leq \frac{C_0}{d(x)}2\min\{h(x),h(y)\}.
\]
Applying the mean value theorem to $t\mapsto\log t$ between $h(x)$ and $h(y)$, we know that there is a number $\xi$ between them such that
\[
 |F(x)-F(y)|
 =\frac{|h(x)-h(y)|}{\xi}
 \leq \frac{|h(x)-h(y)|}{\min\{h(x),h(y)\}}
 \leq \frac{2C_0}{d(x)}.
\]
This proves \eqref{dist} whenever $d(x)$ exceeds an absolute threshold $D$. If $d(x)\leq D$, the global Lipschitz bound gives $|F(x)-F(y)|\leq A\leq A(1+D)/(1+d(x))$, which completes
the proof of \eqref{dist} after adjusting the constant.

Let $R\geq\max\{1,R_u\}$. Notice that both endpoints of every edge $(x,x+e_i)$ with $x\in Q_{2R}$ and $i=1,2$ lie in $Q_{2R+1}$. Moreover, each lattice point is an endpoint of at most four such edges. Thus, by \eqref{weak} at radius $2R+1$, the number of edges with an endpoint outside $G$ is at most
\[
 4|Q_{2R+1}\setminus G|
 \leq \frac{36}{\eta}|\mathcal N\cap Q_{4R+2}|
 \leq C|\mathcal N\cap Q_{6R}|.
\]
Observe that each of these edges contributes at most $A^3$ to \eqref{energy} by global Lipschitz continuity. Additionally, by \eqref{dist}, the edges with both endpoints in $G$ and initial point $x\in Q_{2R}$ satisfying $d(x)>R$ contribute at most
\[
\frac{C|Q_{2R}|}{(1+R)^{3}}\leq \frac CR.
\]
For each $x\in Q_{2R}$ with $d(x)\leq R$, choose a nearest point
$w_x\in\mathcal N$. Then $|x-w_x|_\infty=d(x)$ and $w_x\in Q_{3R}$.
Grouping by $w_x$ and enlarging each inner sum to the whole lattice
gives
\[
 \sum_{\substack{x\in Q_{2R}\\d(x)\le R}}\frac1{(1+d(x))^3}
 \leq \sum_{w\in\mathcal N\cap Q_{3R}}
      \sum_{x\in\mathbb Z^2}\frac1{(1+|x-w|_\infty)^3}.
\]
For each $w$, exactly $(2n+1)^2-(2n-1)^2=8n$ lattice points satisfy $|x-w|_\infty=n$ for each integer $n\geq 1$. Including the term $x=w$, the inner sum therefore equals $1+\sum_{n=1}^{\infty}8n/(1+n)^3<\infty.$ These edges therefore contribute at most $C|\mathcal N\cap Q_{3R}|$. Adding the contributions and using
$Q_{3R}\subseteq Q_{6R}$ proves \eqref{energy}.
\end{proof}

We are now ready to prove Theorem~\ref{thm:main} as a consequence of the previous lemmas that show how the points where $u$ is negative influence the variation of the majorant $F$. The final step is to use this control to bound the growth of $u$ on large squares and make explicit how that bound improves as the density of negative points decreases.
\begin{proof}[Proof of Theorem~\ref{thm:main}]
Let us assume that $\delta\le\delta_U$. If $\mathcal N=\emptyset$, Lemma \ref{lem:lipschitz-energy} shows that $u$ is constant, and both conclusions in the statement follow. Suppose henceforth that $\mathcal N\ne\emptyset$.

Let $R\ge\max\{1,R_u\}$. For each $a \in \mathbb{Z}^2$, we define the function $\widetilde{F}$ in the square $\{a + se_1 + te_2 : s,t \in [0,1]\}$ as
\begin{equation*}
    \widetilde{F}(a + se_1 + te_2) = \begin{cases}
        (1-s-t)F(a) + s F(a+e_1) + tF(a+e_2), &\text{if }s+t \leq 1,\\
        (1-t)F(a+e_1) + (1-s) F(a + e_2) + (s+t-1)F(a+e_1+e_2), &\text{if } s+t\geq 1.
    \end{cases}
\end{equation*}
Notice that the two pieces of this function coincide in the points corresponding to the diagonal $s+t=1$. Moreover, if we consider two adjacent squares, the function also coincides in the shared edge. Therefore, $\widetilde{F}$ is a continuous extension of $F$ to $\mathbb{R}^2$. 

Observe that if $T_1 = \{a + se_1 + te_2 : s,t \in [0,1] \text{ and } s+t\leq 1\}$, then the gradient of $\widetilde{F}$ in the interior of $T_1$ is
\begin{equation*}
    \nabla \widetilde{F} = (F(a+e_1) - F(a) , F(a+e_2) - F(a)).
\end{equation*}
Hence, H\"older inequality yields
\begin{equation*}
\begin{aligned}
    \int_{T_1}|\nabla \widetilde{F}|^3 &\leq \frac{1}{2}(|F(a+e_1) - F(a)|^2 + |F(a+e_2) - F(a)|^2)^{3/2}\\
    &\leq C (|F(a+e_1) - F(a)|^3 + |F(a+e_2) - F(a)|^3).
\end{aligned}
\end{equation*}
Similarly, we have
\begin{equation*}
    \int_{T_2}|\nabla \widetilde{F}|^3\leq C (|F(a+e_1+e_2) - F(a+e_1)|^3 + |F(a+e_1+e_2) - F(a+e_2)|^3).
\end{equation*}
Decomposing each unit square contained in $[-R,R]^3$ into the corresponding two triangles and using \eqref{energy}, we obtain
\begin{equation}\label{eq:EstimateL3norm}
    \int_{[-R,R]^2}|\nabla \widetilde{F}|^3 \leq  C\left(|\mathcal N\cap Q_{6R}|+\frac1R\right).
\end{equation}
For a bounded real-valued function $f$ on a nonempty set $E$, let us define the oscillation of $f$ in $E$ as
\[
\operatorname{osc}_E f
:=\sup_{x\in E}f(x)-\inf_{x\in E}f(x)
=\sup_{x,y\in E}|f(x)-f(y)|.
\]
On the fixed square $[-1, 1]^2$, Morrey's inequality
\cite[\S5.6.2]{EvansPDE}, applied to $v$ minus its mean and
combined with Poincar\'e's inequality, gives
\begin{equation}\label{eq:Morrey}
\operatorname{osc}_{[-1, 1]^2}v\le C\|\nabla v\|_{L^3([-1, 1]^2)}. 
\end{equation}
Let $v(\xi)=\widetilde F(R\xi)$. Since
$\nabla v(\xi)=R\nabla\widetilde F(R\xi)$, a change of variables
gives
\[
 \int_{[-1, 1]^2}|\nabla v(\xi)|^3\, d\xi = R^3 \int_{[-1, 1]^2}|\nabla \widetilde{F}(R\xi)|^3\, d\xi
 =R\int_{[-R, R]^2}|\nabla\widetilde F(x)|^3\, dx.
\]
Taking cube roots in the above equation, using the fact that $\widetilde{F} = F$ in $Q_R$, and applying \eqref{eq:Morrey}-\eqref{eq:EstimateL3norm}, we obtain
\[
 \operatorname{osc}_{Q_R}F
 \le \operatorname{osc}_{[-R, R]^2}\widetilde F
 \le CR^{1/3}
     \left(|\mathcal N\cap Q_{6R}|+\frac{1}{R}\right)^{1/3}.
\]

Since $g\le F$ and $0\in Q_R$, we have
\[
 \log(1+\|u\|_{Q_R})
 =\max_{x\in Q_R}g(x)
 \le F(0)+\operatorname{osc}_{Q_R}F.
\]
Plugging $|\mathcal N\cap Q_{6R}|=(12R+1)^2\delta_{6R}$ and using $(a+b)^{1/3}\le a^{1/3}+b^{1/3}$, we deduce
\begin{align*}
\log(1+\|u\|_{Q_R})&\le F(0)+C+CR^{1/3}(12R+1)^{2/3}\delta_{6R}^{1/3}\\
 &\le F(0)+C+C'R\delta_{6R}^{1/3}.    
\end{align*}
Given that $F(0)$ is independent of $R$, this proves \eqref{eq:sparse-finite} with $C_u=F(0)+C$. Dividing by $R$, taking the upper limit, and using $\limsup_{R\to\infty}\delta_{6R}\le\overline{\delta}^-(u)$, we obtain
\eqref{eq:sparse-rate}.
\end{proof}

\section{Acknowledgements}
J.M. was partially supported by the Simons Foundation Grant $\# 453576$.

\section{AI Disclosure}
The first version of this manuscript \cite{GaitanMontejoGarzonMadrid2026v1} contained only Theorem \ref{thm:polynomial-main}, \cite[Theorem 1.4]{GaitanMontejoGarzonMadrid2026v1}, and Proposition \ref{prop:large-range}. That is, we had given a positive answer to Question \ref{QuestionPositivity} for functions of polynomial growth and for functions of finite translation rank, and we had deduced an exponential upper bound for the growth of discrete harmonic functions that are non-negative on a large portion of $\mathbb{Z}^2$. Based on the proof of \cite[Theorem 1.1]{BuhovskyLogunovMalinnikovaSodin2022}, we thought that there should also be an analogous exponential lower bound for the growth of such functions, which would give a positive answer to Question \ref{QuestionPositivity} in its more general form.
After interacting with ChatGPT-6.2 Astra, we used Theorem \ref{thm:polynomial-main} to reduce the desired exponential lower bound to \eqref{eq:ThreeSquaresIneq} and improved Proposition \ref{prop:large-range} to Theorem \ref{thm:main}. Together, these two bounds allowed us to prove Theorem \ref{thm:QuestionAnswer} and complete the second version of this manuscript one week after the first. We are aware of the preprint \cite{Hu2026v1}, which was submitted to arXiv between the first version and the process of polishing the second version of this manuscript and claims to give a different proof of Theorem \ref{thm:QuestionAnswer} using some ideas from Proposition \ref{prop:large-range}.

The authors acknowledge the use of ChatGPT-6.2 Astra for discussions during the development of the arguments and for assistance in preparing the manuscript. The authors take full responsibility for the final version, including the accuracy and integrity of all mathematical statements and proofs.

\bibliographystyle{acm}

\bibliography{bibl}

\end{document}